\documentclass[11pt]{amsart}
\usepackage{latexsym, mathrsfs} 
\usepackage{calrsfs}
\usepackage{amsfonts,fancyhdr}
\usepackage{amssymb}
\usepackage{wrapfig, graphicx}
\usepackage{amscd}
\usepackage{epsfig}
\newtheorem{theorem}{Theorem}
\newtheorem{lemma}{Lemma}

\newtheorem{kor}{Corollary}
\newtheorem{prop}{Proposition}

\newlength{\figboxwidth}
\newcommand{\lat}{\Z \oplus \Z i}

\newcommand{\abx}{(X,\omega)}
\newcommand{\aby}{(Y,\tau)}
\newcommand{\abas}{(S,\alpha)}

\newcommand{\so}{ {\rm SO}_2(\R) }

\newcommand{\slr}{ {\rm SL}_2(\R) }

\newcommand{\slz}{ {\rm SL}_2(\Z) }
\newcommand{\pslz}{ {\rm PSL}_2(\Z) }
\newcommand{\pslr}{ {\rm PSL}_2(\R) }

\newcommand{\slv}[1]{ {\rm SL}(#1) }
\newcommand{\pslv}[1]{ {\rm PSL}_2(#1) }

\newcommand{\R}{\mathbb{R}}
\newcommand{\C}{\mathbb{C}}
\newcommand{\Q}{\mathbb{Q}}
\newcommand{\Z}{\mathbb{Z}}
\newcommand{\N}{\mathbb{N}}
\newcommand{\pro}{\mathbb{P}}
\newcommand{\proj}{\mathbb{CP}}
\newcommand{\T}{\mathbb{T}}
\renewcommand{\H}{\mathbb{H}}

\newcommand{\sC}{\mathscr{C}}
\newcommand{\sE}{\mathscr{E}}
\newcommand{\sF}{\mathscr{F}}
\newcommand{\sH}{\mathscr{H}}
\newcommand{\sL}{\mathscr{L}}
\newcommand{\sM}{\mathscr{M}}
\newcommand{\sO}{\mathscr{O}}
\newcommand{\sP}{\mathscr{P}}
\newcommand{\sS}{\mathscr{S}}

\newcommand{\sU}{\mathscr{U}}

\newcommand{\CC}{\mathcal{C}}

\newcommand{\CF}{\mathcal{F}}

\newcommand{\CH}{\mathcal{H}}

\newcommand{\CL}{\mathcal{L}}

\newcommand{\CQ}{\mathcal{Q}}

\renewenvironment{itemize}{\begin{list}{$\bullet$} 
{\setlength{\labelwidth}{1cm}
\setlength{\leftmargin}{0.4cm}
\setlength{\rightmargin}{0.3cm}
\setlength{\itemindent}{0mm}
\setlength{\parsep}{0.2ex}
\setlength{\topsep}{0.5ex}
\setlength{\itemsep}{0ex}
\slshape}}{\end{list}}

\def\i{\item}
\newcommand{\inum}[1]{\item[#1]}

\DeclareMathOperator{\area}{area}
\DeclareMathOperator{\aff}{Aff^{+}}
\DeclareMathOperator{\aut}{Aut}
\DeclareMathOperator{\vol}{vol}

\DeclareMathOperator{\D}{D}

\DeclareMathOperator{\id}{id}

\DeclareMathOperator{\per}{Per}
\DeclareMathOperator{\hol}{hol}

\DeclareMathOperator{\pr}{pr}
\DeclareMathOperator{\im}{Im}
\DeclareMathOperator{\re}{Re}

\title[ ]{Moduli spaces of 
branched covers of Veech surfaces I: $\mathbf{d}$-symmetric differentials}   

\address{
 University of Notre Dame,  
 206 Hayes-Healy, 
 Notre Dame, IN 46656}  
\email{schmoll.2@nd.edu}

\author{Martin Schmoll}

\date{\today}

\begin{document}

\begin{abstract}
We give a description of asymptotic quadratic growth rates 
for geodesic segments on covers of Veech surfaces in terms of  
the {\em modular fiber} parameterizing 
coverings of a fixed Veech surface. 
To make the paper self contained we derive the necessary 
asymptotic formulas from the Gutkin-Judge formula.
As an application of the method 
we define and analyze $d$-symmetric elliptic 
differentials and their modular fibers $\sF^{sym}_d$. For given genus $g$, 
$g$-symmetric elliptic differentials (with fixed base lattice) 
provide a $2$-dimensional family of translation surfaces. 
We calculate several asymptotic constants, to establish 
their dependence on the translation geometry of $\sF^{sym}_d$ 
and their sensitivity as $\slz$-orbit invariants. 
\end{abstract}

\maketitle
\tableofcontents


\section{Introduction}

During the last decade major progress 
in calculating asymptotic quadratic growth 
rates of geodesic segments on translation 
surfaces has been made. After 
introducing an asymptotic 
formula for a special type of surfaces 
(see \cite{v2,v3}), 
now called lattice-, or Veech-surfaces,   
Veech found a fairly  
general integral-formula \cite{v4} relating the   
asymptotic constants of a translation surface 
to the $\slr$-orbit closure of the surface in 
the moduli space of translation surfaces. 
To evaluate the so called {\em Siegel-Veech formula}  for a particular translation surface 
one needs to calculate the volume of  
the closure of the surfaces $\slr$-orbit in moduli space, which in general is not known to be 
well-defined, since presently 
there are no general structure 
theorems for $\slr$-orbit classification. 
First assuming the $\slr$-orbit closure of a surface 
is a complex manifold, one can 
rewrite part of the Siegel-Veech constants in 
terms of volume of this orbit closure. To 
complete the evaluation of Siegel-Veech constants 
one has to know the orbit(s) of the 
surface with respect to the action  
(of conjugates) of the unipotent group 
\[
\sU:=\left\{\left[\begin{smallmatrix}1&t \\ 0&1
\end{smallmatrix}\right]:\ t \in \R \right\} \subset \slr.\]   
Specifying an $\sU$-orbit closure classification is known as the {\em Ratner problem} for the moduli space of (translation) surfaces (see \cite{emz}).   
Improvement of Veechs work by Eskin and Masur \cite{em} says that the {\em generic} surface in a stratum admits Siegel-Veech constants. 
Essentially a "generic" surface has the whole stratum as its $\slr$-orbit closure. Strata are complex manifolds and 
certain Siegel-Veech constants for "generic" surfaces  in a given stratum have been calculated in 
\cite{em,emz}.   

This leaves the {\em fundamental 
problem} to tell, if a given surface 
actually belongs to the {\em generic set}. 
Jet another difficulty is   
the {\em complexity} arising in 
volume calculations for "spaces" of 
{\em high genus surfaces}.  
Moduli space volumes rely on counting 
torus coverings and have been studied 
in \cite{emz} using results of \cite{eo}.

A class of examples 
where one understands all $\slr$ (and $\sU$) orbit closures 
are spaces of torus-covers. 
In fact spaces of torus-covers cover  
a homogeneous space themselves,  therefore  
we have a {\em Ratner-Theorem} with 
respect to the action of 
the unipotent group 
$\sU$ on spaces of torus-covers.
Recently Eskin, Marklof and Morris-Witte established  
a Ratner Theorem for spaces of coverings of Veech 
surfaces which are not tori (see \cite{emm}). 

For genus $2$ torus-covers 
generic Siegel-Veech constants 
were calculated in \cite{ems}.  
\medskip\\
So far the literature does not provide   
an applicable method to calculate Siegel-Veech constants  
for surfaces in arbitrary high genus, avoiding 
the above mentioned difficulties 
of the general approach. 
There is motivation to have such 
a method partly related to polygonal billiards, 
but also to study the inner structure of series 
of asymptotic constants as found in \cite{ems}. 
In this paper we develop a method to calculate 
Siegel-Veech constants \cite{s3} for 
moduli spaces of coverings of Veech-surfaces and     
express the constants in terms of the 
translation geometry of the modular fiber 
parameterizing coverings of a fixed Veech surface.  
Finally we define and describe a  
$2$-dimensional (fiber-) space of torus-covers 
in any genus, including an evaluation of 
Siegel-Veech constants for all surfaces 
parameterized by these modular fibers. 
\medskip\\
Recall that a translation surface 
is a Riemann surface $X$ together with an atlas whose 
chart changes are translations in $\C$. Equivalently one can  
consider pairs  $(X, \omega)$, where $X$ 
is a Riemann surface and $\omega \in \Omega(X)$ 
a non-trivial holomorphic $1$-form on X. A translation atlas 
allows to glue together local    
the local $1$-forms $dz$,  obtaining 
a global holomorphic $1$-form $\omega$ on $X$. 
The induced metric has finitely many {\em singular points}, 
or {\em cone points}, which are located in 
$Z(\omega)$, the set where $\omega$ 
vanishes. This allows to consider the set 
$SC_{X}(p_1,p_2)$ of all geodesic segments  
connecting two given cone points $p_1,p_2 \in Z(\omega)$, 
{\em not containing any other cone point}. Goedesic 
segments with this property 
are usually called {\em saddle connections}. 

Closed geodesics on $(X,\omega)$ 
not containing cone points typically 
appear in (maximal) families 
of parallel geodesics having the 
same length, say $w$. 
Any maximal family covers an 
open set on $X$ isometric to a cylinder 
$\CC=\R/w\Z \times (0,h)$, the positive number 
$h$ is called {\em height} of the cylinder and $w=w(\CC)$ 
is its {\em (core) width}. The images of the closed geodesics 
under this isometry are {\em horizontal loops} 
$\R /w\Z \times \{l\} \subset \CC$.
 
To specify the growth rate problems take  
the set of maximal cylinders $Cyl_X$ on $X$   
or one of the sets $SC_X(p_1,p_2)$ and 
consider the counting function  $\R_+ \rightarrow \N$
\[ T \mapsto N(Cyl_X,T):=|\{\CC \in Cyl_X: w(\CC)<T \}|.  \]
The following result is fundamental  (see \cite{m1, m2} 
and \cite{em} for a recent version):
\medskip\\ 
{\bf Theorem} [Masur]
{\it There are constants $c_{min}$, $c_{max}$ 
such that   
\[ 0< c_{min}T^2 < N(Cyl_X,T)<c_{max}T^2, \quad \text{ 
whenever } \ T>>0. \]
} \medskip

\noindent A similar statement holds for the function 
$N(SC_X(p_1,p_2),T)$ counting saddle connections 
between $p_1\in X$ and $p_2 \in X$ shorter than $T$. 
In case the {\em space} $\overline{\slr \cdot (X,\omega)}$  
obeys a Ratner Theorem,  Masur's Theorem 
upgrades to  
\[\pi \cdot c_{cyl}(X)=
\lim_{T \rightarrow \infty}\frac{N(Cyl_X,T)}{T^2}.\]
We recall that the $\slr$-action on translation surfaces 
is defined as follows: 
\medskip\\
Take a natural chart 
$\zeta(p)=\int^{p}_{p_0}\omega$ around $p_0 \in X\backslash 
Z(\omega)$ and simply postcompose this chart 
using $A \in \slr$ acting real linear on $\R^2 \cong \C$. 
Doing this for all natural charts gives a new translation 
structure  
$(\widetilde{X}, \widetilde{\omega}) =A\cdot (X,\omega)$. 
Obviously the numbers and orders of the zeros 
of $\widetilde{\omega}$ agree with the ones of $\omega$. 
Since the $\slr$-action descends to an action on the moduli space 
of translation structures, the above implies that 
the $\slr$-action preserves sets of Abelian differentials 
with given singularity pattern, so called {\em strata}.  
Strata are complex orbifolds \cite{v1} 
%
%
%
%
and $\overline{\slr \cdot \abx}$ is typically (for the generic surface) 
the whole stratum, but sometimes just a subset (of shape which has yet 
to be described). 
\medskip\\
We are interested in Abelian differentials $\abx$ which are 
branched covers of tori.
\medskip\\ 
{\bf Elliptic differentials.}\ A  {\em lattice}  
$\Lambda \subset \R^2 \cong \C$ is a subgroup 
such that $\R^2/\Lambda$ is a torus. 
The pair $\abx$ is called {\em elliptic differential}, 
if there is a lattice $\Lambda$ and holomorphic 
map $\pi: X \rightarrow \C/\Lambda$
, such that
\[\omega = \pi^{\ast}dz.\]
Throughout the paper we assume lattices 
denoted with $\Lambda$ are {\em unimodular}, 
i.e.:
\[ \int_{\C/\Lambda}dx \wedge dy =1.\] 
However for our purpose a slightly extended notion 
of elliptic differentials is appropriate to cover 
{\em degenerated differentials}.
\medskip\\
{\bf Degenerated elliptic differentials.}
An {\em elliptic differential} 
$(X,p_1,...$  $ ,p_n,\omega)$ consists of the following data:
\begin{itemize}
\i a topological space $X$ and a set of points 
$\{p_1,...p_n\}\in X$  
\i a Riemann surface structure on $X-\{p_1,...p_n\}$
\i a branched covering $p: X-\{p_1,...p_n\} \longrightarrow \C/\Lambda$
\i an induced differential $\omega = p^{\ast}dz$
on $X-\{p_1,...p_n\}$.
\end{itemize}
We call an elliptic differential {\em regular},  
if $(X,p_1,...,p_n,\omega)$ extends to a differential 
on $X$, otherwise we call $(X,p_1,...,p_n,\omega)$ 
{\em degenerated}. 
Note that this notion includes for example unions
of tori over the same base lattice, identified in
finitely many points.
This class of degenerated elliptic differentials appears naturally 
in the compactification of modular fibers.
The degenerated differentials are just 
{\em stable curves} equipped with a holomorphic differential 
away from finitely many {\em singular} points. 
\medskip\\
Connected spaces of elliptic differentials $\sE$, which are closed under the 
action of $\slr$ admit a {\em fiber bundle} structure: 
\begin{equation}
\sF  \longrightarrow \sE  \longrightarrow \slr/\slz.
\end{equation}
Here the fiber $\sF$ is the set of all torus covers in $\sE$ 
which cover a 
particular torus $\C/\Lambda$. Usually we consider 
the fiber over $\Lambda = \Z^2$ which causes 
no restrictions of generality.  
It is easy to see, that the {\em lattice of absolute periods} of $\tau$, i.e.  
\[ \per(\tau):=\left\{\int_{\gamma}\tau: [\gamma]\in H_1(Y;\Z) \right\} \subset \C, \] 
equals $\Lambda$ for all $\aby \in \sF$.    
We will prove the following special version  of  a statement in \cite{s2}: 
\begin{theorem}\label{fiber-structure} 
Given an elliptic differential $\aby$ with 
$\per(\tau)=\Z^2$ and exactly 
two cone points $p_0,p_1$.  
Assume  
$\pi(p_1)-\pi(p_0) \notin \Q^2/\Z^2$, i.e. is not 
torsion with respect to the natural 
translation cover $\pi: Y \rightarrow Y/\per(\tau)$. 
Then $\sF_{\tau}:=\overline{\slz\cdot \aby}$  
admits 
a natural structure as an elliptic differential
$(\sF_{\tau}, \omega_{\tau})$
with a lattice group of affine homeomorphisms: 
\begin{equation}
\D\aff (\sF_{\tau}, \omega_{\tau})\cong 
\slv{\C/\per(\tau),dz}\cong \slz. 
\end{equation} 
In addition $\theta \in S^1$ is a periodic direction on $\abx$, if 
and only if it is periodic 
on  $\sF_{\tau}$. 
$\sF_{\tau}$ decomposes into $n_{\sF}$ 
maximal open cylinders $\sC_1,...,\sC_{n_{\sF}}$, 
with top boundary components $\partial^{top} \sC_1,...,
\partial^{top}\sC_{n_{\sF}}$. 
\medskip\\
The cylinder decomposition $\CC_{i,1}(\alpha),...,\CC_{i,n_i}(\alpha)$ 
of any $\abas \in \sC_i$ in direction $\theta$  is the same, in the sense  
that the foliation $\CF_{\theta}\abas$ contains $n_i$ cylinders 
having width  $w_{i,j}=w(\CC_{i,j}(\alpha))$ 
independent of $\abas \in \sC_i$. 
\end{theorem}
\noindent {\bf Remark 1.} The space $\sF_{\tau}$ might not 
be connected,  for this reason we use its   
affine group rather than assigning a 
"Veech group" to it.
For example think of an arithmetic torus 
cover $\aby$, then the modular fiber 
\[ \overline{\slz \cdot \aby}= \sF_{\tau}. \]
is a finite union of Veech surfaces with 
Veech group a sub-group of $\slz$. 
We then say the affine group of this union of surfaces, i.e.  
the modular fiber $\sF_{\tau}$, is $\aff (\sF_{\tau})$ 
which in turn has linear part $\D \aff (\sF_{\tau})$ $=\slz$.  
More general, if we start with a torus cover $\aby$ 
and obtain a $2$-dimensional modular fiber
\[ \overline{\slz \cdot \aby}= \sF_{\tau}= 
\overline{\aff (\sF_{\tau}, \omega_{\tau})\cdot \aby},  \]
having affine 
group  $\aff (\sF_{\tau}, \omega_{\tau})$ 
with linear part $\D \aff ( \sF_{\tau}, \omega_{\tau})$
$\cong \slz$. It is known (see \cite{ghs}, \cite{emm}) that 
the $\aff (\sF_{\tau}, \omega_{\tau}) $-orbit of a non-periodic "point" 
$\aby \in \sF_{\tau}$ is equally distributed 
in $\sF_{\tau}$ (w.r.t. Lebesque measure).\\ 
{\bf Remark 2.}\ Theorem 1 can be generalized 
to non-arithmetic surfaces.  
Modular fibers $(\sF_{\tau}, \omega_{\tau})$ 
are surfaces (or spaces) which should be studied for itself. 
The modular fibers of genus $2$ torus-coverings  
for example provides a family of arithmetic surfaces   
with unbounded genus, but yet accessible structural 
properties, see \cite{s3}. 
Subspaces  of (higher dimensional)  
modular fibers which are Veech-group invariant 
and $2$-dimensional are studied in \cite{hst}. 
  
In this paper we study the easiest case of modular fibers:  
when the modular fiber is a torus itself. 
\medskip\\
{\bf Asymptotic constants via modular fibers.}
From now on, whenever we write down a modular 
fiber as a pair $(\sF,\omega_{\sF})$, i.e. equipped with 
an Abelian differential, we assume it has dimension $2$. 
\begin{theorem}\cite{s3}\label{siegelveech}
Assume $(\sF_{\tau}, \omega_{\tau})$ is
the modular fiber of a space of elliptic differentials 
which are translation covers of $(\R^2/\Z^2,dz)$. 
Denote the maximal horizontal cylinders in   
$\CF_h(\sF_{\tau},\omega_{\tau})$ by  
$\sC_1,...$ $,\sC_{n_{\sF}}$. 
Further assume the $n_i$ maximal cylinders in the 
horizontal foliation of  $\abas \in \sC_i$ 
have width $w_{i,1},...,w_{i,n_i}$. 
If $\abas \in \sF_{\tau}$ has infinite 
$\aff(\sF_{\tau}, \omega_{\tau})$ orbit, 
then it admits the Siegel-Veech constant  
\begin{eqnarray}\label{genericcount}
c_{cyl}(\alpha)=
\frac{1}{\area(\sF_{\tau})}\sum^{n_{\sF}}_{i=1}
\sum^{n_{i}}_{k=1}\frac{\area(\sC_i)}{w^2_{i,k}}
\end{eqnarray} 
for maximal cylinders. If  the orbit
\[\sO_{\alpha}:=\aff(\sF_{\tau}, \omega_{\tau})\cdot \abas 
\subset \sF_{\tau}\] 
 is finite, we have 
\begin{eqnarray}\label{torsioncount}
c_{cyl}(\alpha)=
\frac{1}{|\sO_{\alpha}|}
\sum^{n_{\sF}}_{i=1}\left(  
\sum^{n_{i}}_{k=1}\frac{|\sO_{\alpha} \cap \sC_i|}{w^2_{i,k}}+  \right. 
\left. \sum^{m_i}_{k=1}
\frac{|\sO_{\alpha} \cap \partial^{top} \sC_i|}{w^2_{i,k}}
\right).
\end{eqnarray}
\end{theorem}
The Theorem expresses Siegel-Veech constants 
as quantities depending on the horizontal foliation 
of the modular fiber $(\sF_{\tau}, \omega_{\tau})$.  
Note that each leaf $\sL \in 
\sF_h(\sF_{\tau}, \omega_{\tau})$ has  
a finite $\sU \cap \slz$ orbit. Looking at the 
$\sU$ orbit of $\sL$ in $\sE$, we see either a disjoint 
union of cylinders, if $\sL$ is a saddle connection, or  
a disjoint union of  $2$-tori in $\sE$, if $\sL$ is a 
regular loop.  That shows, the $\sU$ orbit closures 
(of $\abas$),   
might be much smaller than the $\slz$-orbits (of $\abas$).  
\medskip\\
{\bf Saddle connections.} For saddle connections 
one obtains a very similar formula involving the 
saddle connections in $\CF_h(\sF_{\tau},\omega_{\tau})$. 
Indeed given a saddle connection $s$ on a finite 
$\slz$ orbit surface $\abas \in \sF_{\tau}$ 
with two cone points we can always find an $A \in  \slz$, 
such that on $A \cdot \abas \in \sF_{\tau}$ 
the image of $s$, i.e. $A\cdot s$,  is horizontal. 
Restriction to the two cone point situation 
implies, that for each horizontal saddle connection 
$s \in \CF_h  \abas$ which connects one and 
the same cone-point there is a saddle 
connection in 
$s_{\sF} \in \CF_h(\sF_{\tau}, \omega_{\tau})$ 
with $|s|=|s_{\sF}|$.  In particular all 
horizontal saddle connections on 
$(\sF_{\tau}, \omega_{\tau})$ have integer length, 
if $\abas$ covers $\R^2/\Z^2$. 
Moreover each rational 
point on a horizontal saddle connection 
of $\sF_{\tau}$ represents an arithmetic 
surface $\abas$, 
having a horizontal saddle connection $s$ 
which connects the two cone points of $\abas$. 
In fact as a point $\abas \in s_{\sF} \in 
\CF_h(\sF_{\tau}, \omega_{\tau})$ divides 
the saddle connection $s_{\sF}$ into 
two pieces, say $s^{+}$ and $s^{-}$,  
and $\abas$ contains horizontal 
saddle connections of length 
$|s^+|$ and  of length $|s^{-}|$ connecting 
the two different cone points. There are 
two multiplicities $m^{\pm}_{s_{\sF}}$ 
of saddle connections of length $s^{\pm}$ 
which are the same for all $\abas \in s_{\sF} $.   
\begin{theorem}[{\bf Saddle connections}]\cite{s3}
Let $SC_h(\sF_{\tau}) \subset 
\CF_h(\sF_{\tau}, \omega_{\tau})$ be the set of all 
horizontal saddle connections in $\sF_{\tau}$, then 
taking the assumptions and notations of Theorem 
\ref{siegelveech}, 
we find for the asymptotic quadratic growth rate 
$c_{\pm}(\alpha)$ of  saddle connections on $\abas \in \sF_{\tau}$ connecting the two different cone points of $\abas$ 
\begin{eqnarray}\label{torsadcount}
c_{\pm}(\alpha)=
\frac{2}{|\sO_{\alpha}|}\sum_{s \in SC_h(\sF_{\tau})}
 \sum_{(Z,\nu) \in \sO_{\alpha}(s)}\frac{m^+_s}{|s^{+}|^2}.
\end{eqnarray}
with $\sO_{\alpha}(s):=\sO_{\alpha} \cap s$ in the finite orbit case. 
For generic $\abas$ we find 
\begin{eqnarray}\label{gentorsadcount}
c_{\pm}(\alpha)=
\frac{2\zeta(2)}{\area(\sF_{\tau})}\sum_{s \in SC_h(\sF_{\tau})}
m^+_s 
=\frac{2\zeta(2)}{\area(\sF_{\tau})}\sum_{p \in Z(\omega_{\tau})}
m^+_s \cdot \widehat{o}_p,  
\end{eqnarray}
where $\widehat{o}_p=o_p+1$ and $o_p$ is the order 
of the zero $p \in  Z(\omega_{\tau})$.
\end{theorem}
\noindent {\bf Remarks.} One can break down 
these constants into smaller parts associated 
with  the $\sU$ orbit 
of a saddle connection $s_{\sF} \in SC_h(\sF_{\tau},\omega_{\tau})$, 
or $\slz$-orbits of cone points $p \in Z(\omega_{\tau})$ 
in the generic case. 
The terms $m^-$ and $s^-$ do not appear in the expressions for 
the constants,  
because there is an isometric involution 
$\phi \in \aff (\sF_{\tau},\omega_{\tau})$   
of $\sF_{\tau}$ which exchanges "$+$" and "$-$",   
causing the factor $2$ in the expressions.
\medskip\\
Let us also mention, that in case of our example 
the generic 
constant $c_{\pm}(\alpha)$
can be obtained by taking the limit over 
constants of any sequence of arithmetic surfaces 
$(S_i,\alpha_i) \in \sF_{\tau}$, i.e. 
\[\lim_{i \rightarrow \infty} c_{\pm}(\alpha_i)=c_{\pm}(gen)\]
while it is (in general) not true that 
\[\sum_{(Z,\nu) \in \sO_{\alpha_i}(s)}\frac{m^+_s}{|s^{+}|^2} \rightarrow 
\frac{2\zeta(2)}{\area(\sF_{\tau})}m^+_s,  \quad \mbox{when} \ 
i \rightarrow \infty. \]
The reason for that has to do with the 
following observations: 
\begin{itemize}
\item For some saddle connections 
$s_{\sF} \in SC_h(\sF_{\tau}, \omega_{\tau})$ there 
are sequences of surfaces $\{(S_i,\alpha_i)\}^{\infty}_{i=1}$ 
such that $\sO_{\alpha_i} \cap s_{\sF}=\emptyset$ 
for all $i \in \N$. 
\item The $\slz$ orbit of 
all surfaces contained in the set 
\[A := \sF_{\tau}\backslash\  \slz \cdot  
SC_h(\sF_{\tau}) \subset \sF_{\tau}\] 
avoids the set $\CF_h(\sF_{\tau}, \omega_{\tau})$. \\
Note, that  $A$ is obviously an $\slz$-invariant $G_{\delta}$-set. 
\end{itemize}   
\medskip

\noindent{\bf $\mathbf{d}$-symmetric differentials.}\  
The main goal of this paper is to carry out 
a completely explicit example 
illustrating Theorem \ref{fiber-structure}
and evaluating the asymptotic constants in 
Theorem \ref{siegelveech}. 
To begin with,  consider the following 
class of elliptic differentials: 
\medskip\\
We call $\aby$ {\em $d$-symmetric}, if 
\begin{itemize}
\i  $\tau$ has exactly {\em two} zeros of order $d-1$ 
\i  $\Z/d\Z \subset \aut \aby$
\i  $\deg(\pi)=\int_X \pi^{\ast}(dx \wedge dy) = d$\vspace*{1mm} 
\end{itemize} 
with the natural projection $\pi: \aby \rightarrow \C/\per(\tau) $. 
Here $\aut \aby$ is the group 
of affine homeomorphisms of $\aby$, 
preserving the holomorphic one form $\omega$.
The automorphism group is finite, if $g(Y)>1$, 
since homeomorphisms of $Y$ preserving 
$\tau$ are in fact holomorphic.
We obtain $d$-symmetric differentials using 
the following:
\medskip\\
{\bf Connected sum construction.}\ 
Given an Abelian differential 
$\abx$ and a regular leaf $\CL \in \CF_{\theta}(X)$. Take 
$a \in \CL$ and define the line segment 
\[I:=[0,\epsilon]e^{i\theta}+a \subset \CL.\] 
Then for $d \geq 2$ and a cycle $\sigma \in S_d$ we 
define the Abelian differential
\[ \#^d_{I,\sigma} (X,\omega ) \quad \text{or} \quad  
(\#^d_{I,\sigma} X,\#^d_{I,\sigma} \omega )\]
by slicing $d$ named copies $X_1,...,X_d$ of $X$ 
along $I$ and identify opposite sides of the slits 
according to the permutation $\sigma$. 
The notation  $(\#^d_{I,\sigma} X,\#^d_{I,\sigma} \omega )$ 
is useful, if one needs to consider the differential 
$\#^d_{I,\sigma} \omega$, or the Riemann surface 
$\#^d_{I,\sigma} X$. 
\smallskip\\
The differential $\#^d_{I,\sigma} \omega$ on $\#^d_{I,\sigma} X$ 
is uniquely defined by the property 
\[\#^d_{I,\sigma} \omega|_{X_i}=\omega_i=\omega.\] 
Note: we can rename the $d$ copies of $X$ such that 
the cycle $\sigma$ becomes $\tau=(1,2,3,...,d)$. 
In this case we simply write: 
\begin{equation} \#^d_{I} (X,\omega )
= \#^d_{I,\tau} (X,\omega ).
\end{equation}  
We apply this construction to a torus covering, 
by taking $d$ copies of $\T^2$ and cut them along 
the projection of the line segment 
$I=I_v = [0,v]\subset \C$  
to $\T^2$. For the first we assume the image 
of $I_v$ on $\T^2$ is not a loop. 
Denote the result by 
\[ \#^d_{I}(\T^2,dz) \quad \text{or simply} \quad \#^d_{I}\T^2.  
\]
One extends this to all line segments 
$I_v=[0,v] \subset \C$, by choosing  
a sequence $I_{v_i}$ of 
segments projecting injectively to $\T^2$,  
such that $I_{v_i}\rightarrow I_{v}$ when  
$i \rightarrow \infty$ and takes    
\[ \#^d_{I_v}(\T^2,dz)=\lim_{i \rightarrow \infty} 
\#^d_{I_{v_i}}(\T^2,dz)\] 
on Teichm\"uller space.

The torus covers obtained by this connected 
sum construction have genus $d$ and 
the translation structure has precisely two cone 
points of order $d$. 
For $d=2$ we obtain the example from the introduction of 
part one of this paper and for $d=1$ one obtains 
$2$-marked tori studied in 
\cite{s1}.
By construction every $\#^d_{I}(\T^2,dz)$ 
has automorphism group 
\[\Z/d\Z \cong \aut(\#^d_{I}(\T^2,dz)) \subset 
\aff(\#^d_{I}(\T^2,dz))\]
and is a $d$-fold torus cover. 
Therefore all $\#^d_{I}(\T^2,dz)$ are $d$-symmetric.  
\medskip\\  
Denote the {\em set of isomorphy classes of $d$-symmetric   
coverings} of $\T^2=\R^2/\Z^2$ by $\sF^{sym}_d:=\sF^{sym}_d(d-1,d-1)$. Then 
\begin{theorem}\label{symstructure}
The modular fiber $(\sF^{sym}_d, \omega_{d})$ is  
given by 
\[ (\sF^{sym}_d, \omega_{d}) \cong (\R^2/d\Z^2 - 
\Z^2/d\Z^2,dz).\] 
The integer lattice $\Z^2/d\Z^2 \subset \T^2_d:=\R^2/d\Z^2$ 
represents degenerated differentials (with one cone point) 
in the compactification of $\sF^{sym}_d$.  
\smallskip\\
Every $d$-symmetric differential in  $\sF^{sym}_d$
is represented by  $\#^d_{I}(\T^2,dz)$ for an $I \subset \T^2$.
\end{theorem}
\noindent There is an obvious, but interesting consequence: 
\begin{kor}
Given natural numbers $g$ and $n \geq 2$,  there is a $d$-symmetric differential $\aby \in \sF^{sym}_g$, having the congruence 
group $\Gamma_1(n)$ as Veech group. 
In fact an $\aby$ with this Veech group is given by 
\[[\#^d_{g/n}(\T^2,dz)]=[g/n+i0] \in \C/g(\Z\oplus \Z i). \]
\end{kor} 
\noindent Note that $\#^d_{g/n}(\T^2,dz)$ is degenerated, if 
and only if $g/n \in \Z$, or in other words $n|g$. 
Thus given $n\geq 2$, there are always infinitely 
many non-degenerated 
$d$-symmetric differentials with Veech group 
$\Gamma_1(n)$ as well as infinitely many 
degenerated ones with the same Veech group.
\medskip\\
\noindent For $d$-symmetric differentials we evaluate formula 
\ref{genericcount} and find: 
\begin{theorem}\label{siegelveechdsymm}
Let $\abas \in \T^2_d$ be $d$-symmetric and 
$\abas \notin \Q^2/d\Z^2$, i.e. has infinite $\slz$-orbit in $\T^2_d$. 
Then the asymptotic quadratic growth 
rate of periodic cylinders $c_{cyl}$ on $S$ is: 
\begin{equation}
c_{cyl}(\alpha)= 2\sum_{p|d} \frac{\varphi(p)}{p^3}.
\end{equation}
\end{theorem}
\noindent We evaluate formula \ref{torsioncount} for $d$-symmetric 
differentials $\abas \in \Q^2/d\Z^2 $ $\subset \T^2_d$,   
the {\em torsion points} in $\sF^{sym}_d$, as well. Finally we 
analyze special subsets of all saddle connections connecting 
the two cone points on $\abx \in \sF^{sym}_d$.  
It turns out,  that for given $d$ and all 
$1 \leq m \leq d$ 
with $m|d$, there are $a$ chains, each
consisting of $d/m$ parallel saddle connections which 
degenerate simultaneously, 
i.e. have the same length. 
We call such configuration of saddle connections 
$m$-{\em homologous}.
Together all $d/m$ chains belonging to a set 
of $m$-homologous saddle connections 
define a boundary in homology. 
We write down the asymptotic 
constants for $m$-homologous 
saddle connections for given $\abx \in \sF^{sym}_d$. 
\medskip\\
{\bf Convention.} For objects related to modular  fibers  
we will use script fonts like $\sF$, $\sC$ and $\sL$, while 
for (objects on) translation surfaces, which are not modular fibers, 
we use calligraphic fonts $\CC$ and $\CL$. 
\medskip\\
{\bf Acknowledgments.} This paper would not have been 
possible without the encouragement and sponsorship of Serge Troubetzkoy 
(IML-Marseille), Anatole Katok (Penn State) and Marc Burger (ETH). 
The author is especially grateful for the support of the 
Max-Planck-Institute for Mathematics, Bonn in July 2004 and 
for stimulating discussions with Curtis T. McMullen and 
Anton Zorich during this time. 

\section{Background}

{\bf Geodesic segments on elliptic differentials.}
For any $\theta \in [0,2\pi]$ the direction foliations 
$\CF_{\theta}$ given by 
$\ker(\sin(\theta) dx+ i\cos(\theta) dy)$ on $\C$ 
define foliations $\CF_{\theta}(\C/\Lambda)$ on
any torus $\C/\Lambda$. The {\em leaves} of 
$\CF_{\theta}(\C/\Lambda)$ are straight lines and 
geodesics on $\C/\Lambda$ with respect to the 
metric $|dz|^2=dx^2+dy^2$. 
These objects pull back to any elliptic differential 
$\abx$, i.e. there are direction foliations $\CF_{\theta}\abx$ 
having integral curves which are straight lines with respect 
to the metric $|\omega|^2$. 
\medskip\\
Given a covering $\pi: X \rightarrow \C/\Lambda$ with 
$\omega = \pi^{\ast}dz$, then a leaf $\CL \in \CF_{\theta}\abx$ 
is compact if and only if $\pi(\CL) \in \CF_{\theta}(\C/\Lambda,dz)$ 
is compact. 
In particular, if $\abx$ is a cover of $(\T^2,dz):=(\C/\lat,dz)$,  
leaves of $\CF_{\theta}\abx$ are compact 
whenever $\tan(\theta) \in \Q$. 
A compact leaf which either has integral curve $\gamma$, or is represented by a homology cycle $[\gamma]$ is denoted by $\CL_{\gamma}$.    
We call a leaf $\CL$ {\em regular}, if it does not 
contain any zeros of $\omega$.
\medskip\\
{\bf Dynamics of $\slz$ on $\T^2$.} We collect facts on the $\slz$ operation 
on $2$-{\em tori} $\T^2_{[0]}=(\R/\Z^2,dz,[0])$ 
(marked in $[0]$). 
Recall that marked tori are the modular fibers 
for moduli spaces of $d$-symmetric differentials. 

It is clear that the affine group $\aff(\T^2,dz)$ 
is isomorphic to $\slz \ltimes \C/\lat$ with the multiplication 
\[(A,[a])\circ (B,[b])=(A \cdot B, [b+Aa]) \ \mbox{ 
where } a,b \in \R^2  \ \mbox{ and } A,B \in \slz. 
\] 
Clearly the short exact sequence of affine map becomes
\begin{equation}
\begin{CD}
0 @>>> \C/\lat  @>{i}>> \slz \ltimes \C/\lat @>{\D}>> \slz @>>> 1
\end{CD}
\end{equation}
As mentioned above, we mark one point on 
$\T^2$, say the origin $[0]$,  
to get rid of translations $\aut(\T^2,dz) \cong \C/\lat$. 
\medskip\\
By $\slv{\T^2_{[x]},dz} \subset \slz$ we denote the stabilizer of 
the point $[x]\in \T^2$. The projection of $\slv{\T^2_{[x]},dz}$ 
to $\pslr$ is the Veech group of the two marked torus 
$\T^2_{[x]}=(\R/\Z^2,dz,[0], [x])$ and we can regard  
$\slv{\T^2_{[x]},dz}$ as the Veech group of $(\T^2_{[x]},dz)$  
with distinguished marked points. 
If and only if $[x]=x+\Z^2$ is rational 
$\slv{\T^2_{[x]},dz}$ is a lattice in $\slz$, in fact:
\begin{prop}
Given $a,b,n \in \N_0$ with $\gcd(a,b,n)=1$ then the stabilizer  
$\slv{\T^2_{[x]},dz}$ of $[x]=[\frac{a}{n},\frac{b}{n}]$ 
is conjugated to 
\[ \Gamma_1(n):=\left\{
\left( 
\begin{smallmatrix} 
 a & b \\ c & d
\end{smallmatrix}\right)
\in \slz: \   
\left( 
\begin{smallmatrix} 
 a & b \\ c & d
\end{smallmatrix}\right)
\equiv \left( 
\begin{smallmatrix} 
 1 & b \\ 0 & 1
\end{smallmatrix} \right)
\! \mod \ n \right\} 
\subset \slz.\]
Moreover for $[x]=[\frac{1}{n},0]$: 
$\slv{\T^2_{[x]},dz}=\Gamma_1(n)$. 
\end{prop} 
\noindent This is easy to see, compare \cite{s1} up to the fact that 
there we allow affine maps which interchange the marked points. 
The linear part of this bigger affine group is  
isomorphic to $\Gamma^{\pm}_1(n):=\{\pm A: A \in \Gamma_1(n)\}$. 
\medskip\\
The $\slz$-orbit of $[\frac{1}{n},0] \in \T^2$ is the set 
\begin{equation}\label{torbit} 
\begin{split}
\sO_n =\left\{\left[\frac{a}{n}+i\frac{b}{n}\right] \in \T^2:\  
a,b,n \in \Z \mbox{ with } \gcd(a,b,n)=1 \right\}=\\
=\slz \cdot \left[\frac{1}{n}\right],
\end{split}
\end{equation} 
in particular 
\[ 
[\Gamma_1(n):\slz]=|\sO_n|= 
n^2 \prod_{p|n}
\left(1-\frac{1}{p^2} \right)=
\varphi(n)\psi(n).\] 
The product is taken over all prime divisors $p$ of $n$.  
We have also used the {\em Euler $\varphi$ function} 
and the {\em Dedekind $\psi$ function}:
\begin{equation}
\varphi(n):=n \prod_{p|n}
\left(1-\frac{1}{p} \right), \quad 
\psi(n):= n \prod_{p|n}
\left(1+\frac{1}{p} \right).
\end{equation}
\vspace*{2mm}\\
{\bf Torsion points on ${\mathbf \T^2_d}$.}
A point on $\T^2_d$, contained in  
\[ \T^2_d[m]:=\ker(\T^2_d\stackrel{m}{\longrightarrow}\T^2_d)=
\frac{d}{m} \Z^2/d(\Z^2).\]
for some $m \in \N$ is called {\em torsion} point. 
The homomorphism denoted by $m$ is 
$z \mapsto m\cdot z$ for $z \in \T^2_d$.  
The {\em order} of a torsion point $z \in \T^2_d$ 
is the smallest $m\in \N$ such that $m \cdot z = 0 \in \T^2_d$. 
We denote the set of torsion points of order $m$ on 
$\T^2_d$ by $\T^2_d(m)$.  
On the standard torus $\T^2=\T^2_1$ we have by equation \ref{torbit} 
\begin{equation}\T^2(m)=\slz \cdot [1/m]=\sO_m. \end{equation} 
Now the isomorphism 
\[\T^2_d\stackrel{d^{-1}}{\longrightarrow} \T^2, \quad 
z \mapsto d^{-1}\cdot z \] 
is $\slz$ equivariant and 
identifies torsion points of order $m$. 
\medskip\\
{\bf Cusps and fundamental directions.} 
Take a lattice surface $\abx$, then Veech group $\slv{X,\omega}$ 
acts on the set of periodic directions $\pro^1(\hol\abx)$ of $\abx$.  
A representative $v_i \in \pro^1(\hol\abx)$  
for the orbit $\pro^1(\slv{X,\omega}\cdot v_i)$ is called a {\em 
fundamental direction}. The set of fundamental directions 
is finite and in bijection to the cusps of $\slv{X,\omega}$.


\section{Asymptotic formulas for branched covers of Veech surfaces}

Now we develop a formalism for evaluating quadratic 
growth rates for branched covers of Veech surfaces. 
In case the covering is itself Veech, we derive the asymptotic formula   
based on the formula found by Gutkin and Judge \cite{gj}
and independently by Vorobetz \cite{vrb}. In \cite{s2} 
we  use the Siegel-Veech formula to obtain Siegel-Veech 
constants for any dimensions of $\sF$.   
\medskip\\
We start with the following data: 
\begin{flushleft}
\begin{itemize}
\inum{1.} a Veech-surface $(X,\omega)$ with Veech-group $\slv{X,\omega}$ 
\inum{2.} a branched (translation) cover 
$p: \aby \rightarrow \abx$ with lattice 
Veech group \[\slv{Y,\tau} \subseteq \slv{X,\omega}.\]
\inum{3.} a space of covers $\sF_{\omega}$ of $\abx$, such that 
 $\overline{\slv{X,\omega}\cdot \aby} \subseteq \sF_{\omega}$
\end{itemize} 
\end{flushleft}
{\em Note:} in general assumption $2$ is a condition 
(which can always be reached by naming some objects on the cover). 
For arithmetic differentials $\aby$ however, the 
(absolute) period lattice $\Lambda=\per(\tau)\in \C$ 
gives the covering map $p: \aby \rightarrow (\C/\Lambda,dz)$, 
implying $\slv{Y,\tau} \subset \slv{\C/\Lambda,dz}$.  

For an arbitrary covering $p: \aby \rightarrow \abx$ 
we can achieve 
\[\aff\aby \cong \slv{Y,\tau} \] 
by taking the reduced translation structure, 
\[ (Y_{red},\tau_{red}):=(Y/\aut(Y,\tau), \pi_{\ast}\tau)\]
where $\pi:Y \rightarrow Y_{red}$ is the natural   projection.    
\smallskip\\
We continue with assumptions;
\begin{itemize}
\inum{4.} let $v_0,...,v_n \in \pro^1(\slv{X,\omega})$ 
be a {\em set of fundamental directions} 
on the base surface $\abx$. 
Without restrictions we assume $v_0$ is horizontal. 
\inum{5.} furthermore let $C_{i,1},...,C_{i,k_i}$  
are the maximal cylinders in direction $v_i$ on $\abx$.  
\end{itemize}
\vspace*{2mm}

\noindent The following asymptotic formula has been derived by Gutkin and Judge \cite{gj}, compare also \cite{emm}.
\begin{prop}
Assume the stabilizer $N \subset \slr$ of $v \in \R^2$ 
has nontrivial intersection with a lattice $\Gamma \subset \slr$ 
and let $A$ be a generator of $N \cap \Gamma$, then 
\begin{equation}\label{count} \mid \Gamma \cdot v \cap B(T)\mid \sim 
\vol(\H/\Gamma)^{-1}
\frac{<A v^{\perp},v>}{|v^{\perp}||v|^3} T^2, 
\end{equation} 
where $v^{\perp}$ is perpendicular to $v$. 
\end{prop}
\noindent Without restrictions we take
 $v^{\perp}=\left[\begin{smallmatrix}0 \\  h\end{smallmatrix}\right]$ and  
 $v =\left[\begin{smallmatrix}w \\  0\end{smallmatrix}\right] \in \R^2$
and might think of it as a cylinder of height $h$ and width $w$ 
contained in the horizontal foliation of $\abx$. 
Now formula \ref{count} reads 
\begin{equation}\label{count2} \mid \slv{X,\omega} \cdot v\ \cap \ B(T)\mid \sim 
\vol(\H/\Gamma)^{-1}
\frac{<u_{v} v^{\perp},v>}{h \cdot w^3} T^2=
\vol(\H/\Gamma)^{-1}\frac{l_v}{w^2}T^2, 
\end{equation} 
with a linear Dehn twist given by  
\[ u_v= \left[\begin{smallmatrix}1& l_v \\ 0 & 1\end{smallmatrix}\right] \in 
\slv{X,\omega}\]
where $l_v\frac{h}{w}=p/q \in \Q$.   
\medskip\\
Now we relate the above formula for the base surface $\abx$ 
to the one for the cover $\aby$. We already saw: 
\begin{itemize}
\i There is an operation of $\slv{X,\omega}$ on the set of coverings 
$\aby \rightarrow \abx$. In particular there is a transitive operation on the orbit $\slv{X,\omega}\cdot \aby$ in moduli space $\sF_{\tau}$.
\i a direction $v$ on $\abx$ is completely periodic,  
if and only if it is completely periodic on $\aby$ 
\i for a completely periodic direction $v$ the 
stabilizer $N_v\aby \subset \slv{Y,\tau}$ on $\aby$ is a 
subgroup of the stabilizer $N_v\abx \subset \slv{X,\omega}$ and  
\begin{equation}
[N_v\abx:N_v\aby]= \mbox{ length of the orbit } N_v\abx \cdot \aby=:i_v(Y,X).
\end{equation}  
\end{itemize}
The following observation is fundamental for all 
calculations and formulas in this section.
\medskip\\
{\bf Combinatorial description of cusps on $(Y,\omega)$.}
Take the base Veech surface $(X,\omega)$ and choose 
a fundamental direction $v$, again we  
{\em assume $v$ is horizontal}. The fundamental 
directions $v_0,...,v_n$ of a Veech surface $(X,\omega)$  
represent the cusps of $\H/\slv{X,\omega}$. 
\begin{prop}\label{anton}
The number of cusps of $(Y,\nu)$ above a given  
cusp $v$ on $(X,\omega)$ are in one-to-one 
correspondence to the $u_v$ orbits on the set 
$\slv{X,\omega}/\slv{Y,\nu}$.    
\end{prop}
\begin{proof}
For every direction $w \in \pro^1(\slv{X,\omega}\cdot v)$ there 
is (at least one) surface 
\[(S,\alpha)= A\cdot (Y,\nu) \in \slv{X,\omega}\cdot(Y, \nu)\] 
for which $A\cdot w$ is horizontal. 
If there are two maps $A,B \in \slv{X,\omega}$ making 
the direction $w$ horizontal, then 
$A\circ B^{-1} \in N_v\abx$. Thus the number of 
all $A\cdot (Y,\nu)\in \slv{X,\omega}\cdot(Y, \nu)$ such that 
$A\cdot w$ is horizontal equals $i_v(S,X)$, which is the length of 
the $u_v$ orbit of $(S,\alpha)$. This proves the claim.
\end{proof}
\noindent We define the {\em relative width of a  
cusp} represented by a direction $w$ 
over $v$ on $(Y,\nu)$ as the length 
of the $u_v$ (the stabilizer of $v$ in 
$\slv{X,\omega}$) orbit 
\[ \{u^n_v \cdot (Y,\nu):\ n \in \Z\} \subset 
\slv{X,\omega}\cdot(Y, \nu). 
\medskip\\  
\]
{\em Note:} for arithmetic surfaces (= torus covers) the relative 
cusp width is the same as the usual cusp width. 
Also for arithmetic surfaces 
the above characterization of cusps was  
discovered by A. Zorich and published in \cite{hl}. 
\medskip\\
Together with equation \ref{count2} we have:
\begin{prop}
Under the above assumptions,let $v \in \C$ be a 
periodic direction on $\abx$ 
and $(S,\alpha)$ a differential in 
the orbit $\sO_{\tau}:=\slv{X,\omega} \cdot (Y,\tau)$, 
then 
\begin{equation}\label{rate1}
\mid \slv{S,\alpha}\cdot v \cap B(T)\mid \sim 
\vol(\H/\slv{X,\omega})^{-1}
\frac{l_v \cdot i_v(S,X)}{[\slv{X,\omega}:\slv{Y,\tau}]}
\frac{1}{|v_i|^2}
\end{equation}
\end{prop}  
Again assume $v$ is a periodic direction on $\abx$,  
pick a periodic cylinder $c_X \subset \CF_v\abx$, 
(or saddle connection $s_X \in \CF_v \abx$).
Take a pre-image $c_S$ (or $s_S$) 
of $c_X$ ($s_X$ respectively) on a differential 
$(S,\alpha) \in \sO_{\tau}$, then by \ref{rate1} 
we have the asymptotic quadratic growth rate 
\begin{equation}\label{rate2}
c_{cyl}(c_S)=\frac{\pi}{3\vol(\H/\slv{X,\omega})} 
\frac{l_v \cdot i_v(S,X)}{[\slv{X,\omega}:\slv{Y,\tau}]}
\frac{1}{|w_{c_S}|^2}.
\end{equation}
using $\sU_{\alpha}:=N_v\abx\cdot \abas$ and 
$|\sO_{\tau}|=[\slv{X,\omega}:\slv{Y,\tau}]$ we get
\begin{equation}\label{rate3}
c_{cyl}(c_S)= \frac{\pi}{3\vol(\H/\slv{X,\omega})} 
\frac{l_v}{|\sO_{\tau}|}
\sum_{(S,\alpha) \in \sU_{\alpha}} \frac{1}{w^2}
\end{equation} 
Now for any fundamental direction $w$ 
of $\aby$ such that $w \in \slv{X,\omega}\cdot v$ 
there is a differential $\abas =A \cdot \aby \in \sO_{\tau}$ 
such that $v=A\cdot w$. Thus we can apply equation \ref{rate3} 
to all cylinders in a cylinder decomposition in direction $w$ 
on $\aby$ by looking at the cylinder decomposition $\CC_{\alpha,1},...,
\CC_{\alpha,n_{\alpha}}$ 
in direction $v=Aw$ on $\abas = A \cdot \aby$. Denote the 
width of the cylinder $\CC_{\alpha,i}$ by $w_{\alpha,i}$, 
then the  {\em asymptotic quadratic growth rate of all 
cylinders on $\aby$ having direction 
contained in the set} 
$\slv{X,\omega}\cdot v$ is given by:   
\begin{equation}\label{rate4}
 c_{cyl,\tau}(v)=\frac{\pi}{3\vol(\H/\slv{X,\omega})}  
\frac{l_v}{|\sO_{\tau}|}
\sum_{\abas \in \sO_{\tau}}
\sum^{n_{\alpha}}_{i =1} \frac{1}{w^2_{\alpha,i}}.
\end{equation}  
\vspace*{1mm}

\noindent{\bf Siegel-Veech constants -- finite orbit case.} 
The Siegel-Veech constant 
for cylinders (or saddle connections) on any 
$\abas \in \sO_{\tau}$ 
is now obtained by taking the sum over all cylinders  
a fundamental directions $v_1,...,v_n$ 
on $(X,\omega)$: 
\begin{equation}\label{rate5}
c_{cyl,\tau}=\frac{\pi}{3\vol(\H/\slv{X,\omega})}  
\frac{1}{|\sO_{\tau}|}
\sum_{\abas \in \sO_{\tau}}\sum^{n}_{k=1}l_{v_k}
\sum^{n_{\alpha,k}}_{i =1} \frac{1}{w^2_{\alpha,k,i}}
\end{equation}  
Note that away from the {\em hyperbolic volume, an index and 
length of cylinders in a particular direction} one needs to 
know only the {\em set of fundamental directions on the base}   
$(X,\omega)$ and the number $l_{v_i}$ associated to 
each fundamental direction.
\medskip\\  
{\bf Basic examples -- arithmetic surfaces.}
Arithmetic surfaces $\aby$ possess a natural map 
to a torus $\C/\Lambda$. Taking an $\slr$ deformation of $\aby$
we might assume $\Lambda=\lat$, the deformation does not change 
quadratic growth rates.  Now on $\T^2=\C/\lat$ there 
is only one fundamental direction, we can assume 
this direction is horizontal and $l_{v_i}=1$.  
We write down the statements for several sub-types of 
saddle connections later, there are many interesting cases 
labeled by various properties to consider.
\medskip\\  
Note that for arithmetic surfaces the formula is completely 
general. We apply it in this paper after 
using the  cylinder decomposition of the moduli 
space fiber to simplify and organize information.   
\medskip\\  
{\bf Basic examples -- non arithmetic Veech surfaces.}
Non arithmetic Veech surfaces cannot have a direction 
with one single closed cylinder. Thus the best situation 
is $(X,\omega)$ has
\begin{itemize} 
\i only {\em one fundamental direction} $v$ with, say $k$ cylinders and  
\i the moduli $w_1/h_1=...=w_k/h_k=:m$ of the cylinders are all the same 
\end{itemize}
here $w_i$ is the width and $h_i$ the height of the i-th cylinder. 
It follows $l_v=m$ and the the Siegel-Veech 
constant for cylinders is
\begin{equation}\label{nonarith} c_{cyl,\omega}= 
\frac{\pi \cdot l_v}{3\vol(\H/\slv{X,\omega})}
\sum^{k}_{i=1}\frac{1}{w^2_{i}}= 
\frac{\pi}{3\vol(\H/\slv{X,\omega})} 
\sum^{k}_{i=1}\frac{1}{h_iw_{i}}.
\end{equation} 
\medskip\\  
A new paper by Eskin, Marklof and Morris \cite{emm} contains 
a discussion of the Veech surfaces $X_n$ obtained by linearization 
of the billiard in the triangle with angles 
$(\pi/n,\pi/n,(n-2)\pi/n)$. A result of Veech \cite{v3}
shows $X_n$ has only one fundamental direction $v$ 
and the cylinder moduli are all the same. 
Thus formula \ref{nonarith} applies (compare \cite{emm} pgs 26-28). 
In fact one can simplify using explicit 
expressions for the $w_i$ and $h_i$, see \cite{v3, emm}.
\medskip\\
Both classes of examples are prototypes 
to consider branched covers above them. 
For a very particular case, $2$-fold covers of $X_n$, this 
has been done in \cite{emm} to get information 
on some non Veech triangular billiards. 
\medskip\\ 

\section{Counting using the horizontal foliation of 
the modular fiber.}

One can organize the asymptotic formula \ref{rate5} 
using modular fibers. Assume $\abx$ is Veech with 
lattice group $\slv{X,\omega}$ and $\pi: \aby \rightarrow \abx$ 
is a (branched) covering.  
Define 
\[\sE_{\tau}:=\overline{\slr \cdot \aby}\subset \sH_g
\mbox{ and }\ \sF_{\tau}:=
\overline{\slv{X,\omega} \cdot \aby}.\] 
The {\em closure} of these two spaces 
is taken within the set of Abelian differentials having 
the same zero order pattern as $\tau$. 
There is another closure, denoted by $\sF^c_{\tau}$ 
($\sE^c_{\tau}$ respectively),  
coming from the compactification of the fibers 
\[A \cdot (\sF_{\tau}, \omega_{\tau}), \quad A \in \slr\] 
viewed as Abelian differentials. 
The set $\sF^c_{\tau}-\sF_{\tau}$ 
parameterizes {\em degenerated differentials}, 
but eventually {\em not} up to isomorphy.
In our definition of $\sE_{\tau}=\sE_{\tau}(o_1,...,o_n)$ 
we distinguish cone points, which makes $\sF_{\tau}$ 
into an Abelian 
differential, if $\dim (\sF_{\tau})=2$ (see \cite{s2}). 
In case $\aby$ is Veech itself 
\begin{equation}\sE^c_{\tau}=\sE_{\tau}\cong \slr/\slv{Y,\tau}
\subset \Omega\sM_g \end{equation} 
and by taking base points with respect to 
$\pi: \Omega\sM_g \rightarrow \sM_g$ 
we obtain a Teichm\"uller curve in $\sM_g$. 
\medskip\\
In this paper we are mainly interested in cases, when 
the dimension of $\sF_{\tau}$ equals $0$ or $2$. 
Since we understand the $0$-dimensional case already, 
we develop the counting theory for $2$-dimensional fibers. 
\medskip\\
\begin{proof}[\bf Proof of Theorem \ref{fiber-structure}] 
For every $\abas \in \sF_{\tau}$ there is a  
map 
\[\pi: \abas \rightarrow \abx\]
which for arithmetic $S$ becomes
\[\pi: \abas \rightarrow (\C/\per(\alpha),dz).\] 
Now choose a symplectic base of the relative homology  
$H_1(Y,Z(\alpha);\Z)$, 
say $(\gamma_1,...,\gamma_{2g}, \gamma_{2g+1},...,\gamma_{2g+n})$.  
Given an ordering $z_0,...,z_n$ of the zeros $z_i \in Z(\alpha)$ 
the {\em relative cycles} $(\gamma_{2g+1},...,\gamma_{2g+n})$ 
can be chosen to have the property  
\[\partial \gamma_{2g+i}=z_0 - z_i.  \]
To use the given symplectic base for other surfaces, we 
suppose every differential $\abas$ comes with a marking 
$\phi_S: (Y, Z(\tau)) \rightarrow (S, Z(\alpha))$, i.e. 
an orientation preserving homeomorphism mapping 
$Z(\tau)$ bijectively onto $Z(\alpha)$.  
This allows to identify $H_1(Y,Z(\tau);\Z)$ with 
$H_1(S,Z(\alpha);\Z)$.   
Using the (relative) homology base,  
one defines local coordinates of  
the connected component of the stratum 
$\sH_1(o_0,...,o_n)$ containing $\abas$ 
\begin{equation}\label{coordinate} 
\abas \mapsto  
\left(\int_{\gamma_1}\alpha,...,\int_{\gamma_{2g}}\alpha,  
 \int_{\gamma_{2g+1}} \alpha,...,\int_{\gamma_{2g+n}} \alpha \right) \subset \C^{2g+n}, 
 \end{equation}
giving the stratum a complex manifold structure.  
The $\slr$ action on differentials 
\[ A \cdot \alpha:=
\left( 
\begin{smallmatrix} 
 1 \\ i
\end{smallmatrix} \right)
\left( 
\begin{smallmatrix} 
 a & b \\ c & d
\end{smallmatrix}\right)
\left( 
\begin{smallmatrix} 
 \im \alpha \\ \re \alpha
\end{smallmatrix}\right)
 \]
obviously commutes with integration 
\begin{equation}\label{commuteint} 
\int_{\gamma}A \cdot \alpha =A \cdot 
\int_{\gamma}\alpha, 
\end{equation}
giving the standard real linear action 
of $\slr$ on coordinates.
\medskip\\
If $\abas$ is arithmetic, the 
absolute periods $\gamma \in H_1(Y;\Z)$ 
generate a lattice 
\[\per(\alpha):=\left\{\int_{\gamma}\alpha: 
\gamma \in H_1(Y,\Z)\right\} 
\subset \R^2\]
and therefore a map of differentials 
\[ (S,\alpha) \rightarrow (\C/\per(\alpha),dz).\]
By definition of the $\slr$-action on differentials 
\[\per(\slv{\C/\per(\alpha),dz} \cdot \alpha)=\per(\alpha),\]
which implies that $\slv{\C/\per(\alpha),dz}$ acts on 
the set of differentials $\abx \in $ $\overline{\slr \cdot \aby}$ 
with absolute period lattice 
$\per(\omega)=\per(\alpha)$.   
For simplicity we assume from now on $\per(\alpha)\cong \Z^2$ and thus all differentials contained in the modular fiber
\[\sF_{\tau}:=\overline{\slz\cdot(Y,\tau)}\] 
have absolute period lattice $\Z^2$, by continuity of local  
coordinates (\ref{coordinate}). 
The difference of  any two {\em relative periods} 
$\gamma_{1}, \gamma_{2} \in H_1(Y,Z(\tau);\Z)$ 
with 
\[\partial \gamma_{1}=z_0 - z_i=\partial \gamma_{2}\] 
is an absolute period, hence there is a well-defined map 
\begin{equation} 
\pi: \left\{ 
\begin{split}
\sF_{\tau} &\longrightarrow  \hspace{1.7cm} \T^2 \times \cdots \times \T^2 \\
\abas &\longmapsto (\int^{z_1}_{z_0}\alpha \hspace{-.3cm} \mod \Z^2 , 
..., \int^{z_n}_{z_0}\alpha \hspace{-.3cm}\mod \Z^2)
\end{split}
\right.
\end{equation}
whenever $\gamma_i$ is a relative cycle with $\partial \gamma_{i}=z_0 - z_i $. 
By definition of local coordinates this map is holomorphic and by construction of $\sF_{\tau}$ 
$\slz$-equivariant. In fact, since the absolute period lattice 
is the same for all surfaces parameterized by  $\sF_{\tau}$, 
the absolute period coordinates are locally constant on 
$\sF_{\tau}$ and thus $\pi$ provides local coordinates 
of $\sF_{\tau}$. 
This implies together with the $\slz$ equivariance  that $\pi: \sF_{\tau} \rightarrow  \pi(\sF_{\tau}) 
\subset (\T^2)^n $ 
is a covering with image a union of tori 
of dimension equal of $\dim(\sF_{\tau})$. 
\medskip\\
By assumption we are interested in 
differentials $\aby$ having exactly $2$ cone points 
covering $\T^2$ in a way that the relative location 
of the cone points is torsion. 
This means we obtain a surjective cover 
\[\pi: \sF_{\tau} \rightarrow  \T^2,\]
since the $\pi$ image of the relative period coordinate is not closed under $\slz$. 
Thus $\sF_{\tau}$ is an arithmetic surface, 
by construction it is square-tiled. 
  
Now we can pull back 
$dz$ to get a holomorphic $1$-form $\omega_{\tau}=\pi^{\ast}dz$ 
on $\sF_{\tau}$, thus obtaining an elliptic differential 
\[(\sF_{\tau}, \omega_{\tau}) \quad \text{with }\quad \per(\omega_{\tau}) \subset \Z^2, \]
which has by definition $\D\aff(\sF_{\tau}, \omega_{\tau})=\slz$. 
By means of its definition as a stabilizer group 
one would not call $\slz$ the Veech group of  
$(\sF_{\tau}, \omega_{\tau}) $,  
if $\sF_{\tau}$ is not connected.  
On the other hand each deformation of a surface $\abas \in \sF_{\tau}$ 
along a path $t \mapsto \gamma(t) \in \slr$,  
such that $\gamma(0)=\id$ 
and $\gamma(1)=A \in \slz$ 
represents a homotopy class 
\[[\gamma] \in \pi_1(\slr/\slz) \cong \slz,\] 
which in turn defines an affine homeomorphism 
\begin{equation}
\phi_{\gamma}: \left\{\begin{split}
\sF_{\tau} &\rightarrow \sF_{\tau}\\
\abas &\mapsto A\cdot \abas.
\end{split}
\right. 
\end{equation} 
That the map is affine with derivative $A \in \slz$, follows immediately 
from the previous definition of the translation structure on $\sF_{\tau}$ and equation 
(\ref{commuteint}). 
Thus we get isomorphisms
\[ \D \aff(\sF_{\tau}, \omega_{\tau}) \cong \slz \cong \pi_1(\slr/\slz).\] 

\noindent {\bf Cylinder decomposition of $\mathbf{(\sF_{\tau}, \omega_{\tau})}$.}\  
As an cover of $\T^2$ the differential 
$(\sF_{\tau}, \omega_{\tau})$ is completely 
periodic in each {\em rational} direction, 
in particular in the horizontal direction. 
Thus we only need to show, that 
every $\abas \in \sF_{\tau}$ admits a completely periodic horizontal 
foliation. But this is obvious, because  
\[\per(\alpha) = \Z^2 \supset \per (\omega_{\tau}).\]
The above generalizes in an obvious way to any deformed fiber 
$A\cdot (\sF_{\tau}, \omega_{\tau})$,  $[A] \in \slr/\slz$. 
Moreover one can extend this result to $2$-dimensional, $\slv{X,\omega}$-invariant subspaces in $X^n$, for any (reduced) Veech surface $\abx$ 
(see \cite{hst}). 
\medskip\\
{\bf Stability of cylinder-decompositions.}
The case we are interested in is $\aby$ 
is a non-arithmetic cover $\pi: Y \rightarrow \T^2$ 
with exactly $2$ cone points $p_0$ and $p_1$,  
we assume $\pi(p_0)=[0] \in \T^2$. 
This implies 
the modular fiber $\sF_{\tau}:=\overline{\slz \cdot \aby }$
has dimension $2$.  
A relative period, say $\gamma_1 \in 
H_1(Y, \{p_0,p_1\};\Z)$, with $\partial \gamma_1=p_1-p_0$,  
provides local coordinates 
\[\abas \mapsto z(p_1)=\int_{\gamma_1}\alpha\]  
 on $\sF_{\tau}$. 
If we choose any absolute cycle $[\gamma] \in 
H_1(Y;\Z) \subset H_1(Y, \{p_0,p_1\};\Z)$, 
the function 
\begin{equation}\label{locconst} \sF_{\tau} \ni \abas  
\mapsto \int_{\gamma}|\alpha| \in \N 
\end{equation} 
is continuous, hence constant on connected 
components of $\sF_{\tau}$. In particular, 
if $\gamma$ is represented by a {\em regular leaf} 
of $\CL \in \CF_{\theta}\abas$, i.e.
\begin{equation}\label{regleaf}   
\int_{\CL}|\alpha|= |\int_{\CL}\alpha| \in \N 
\end{equation}
there is a 
neighborhood of $\abas \in \sF_{\tau}$  
where $\gamma$ is represented by a 
regular leaf of the same length as $\CL$.  
Indeed regular leaves 
$\CL \in \CF_{\theta}\abas$ 
have positive distances from the cone point $p_1$,  
and all these distances are continuous functions 
depending on (the coordinates of) 
$\abas \in \sF_{\tau}$. 
Since the length function \ref{locconst}  
is constant on each connected component 
of $\sF_{\tau}$ for  
curves $\gamma$ such that 
$[\gamma] \in H_1(Y;\Z)$, the only possible 
thing which can happen under deformation 
is that $\gamma$ looses the property  
of being realized by a regular leaf 
(if it is at all). This in turn happens only if 
the cone point $p_1$ crosses the leaf $\CL$ realizing $\gamma$.  If we now fix a 
periodic direction $\theta \in S^1$ on 
$\abas \in \sF_{\tau}$, the same direction 
is is periodic on $\sF_{\tau}$. Now assume 
$(S_t,\alpha_t) \in \sL \in \CF_{\theta}(\sF_{\tau}, \omega_{\tau})$,  
$\sL$ regular and $\{(S_t,\alpha_t): \ t \in \R \} = \sL$. 
Then by the above all regular 
leaves of any $(S_t,\alpha_t) $ in direction 
$\theta \in  S^1$ stay regular since we 
move the cone point $p_1$ in direction $\theta$. 
This completes the proof.
\end{proof}

\noindent{\bf Elliptic differentials.}\ As in the previous proof, we assume 
$\sF_{\tau}=\sF_{\tau}(o_1,o_2)$ is a fiber space 
of {\em elliptic differentials} $\abas$ 
with two cone points of ord er $o_i$ 
and period lattice $\per(\alpha) = \Z^2$.  

Each horizontal cylinder $\sC_i$ of $\sF_{\tau}$ 
is bounded by a union of saddle connections $\partial^{top}\sC_i$ 
parametrizing elliptic differentials 
$\abas \in \sF_{\tau}\cap \partial^{top}\sC_i$ where the two 
cone points $z_0,z_1 \in Z(\alpha)$ are connected by (at least one) 
horizontal saddle connection, say $s \in \CF_h\abas$. 
Thus deforming $\abas$ along $\partial^{top}\sC_i$ 
is nothing but collapsing the two cone points of $\abas$ into one. 
This means points on $\partial^{top}\sC_i$  contain (degenerated) 
elliptic differentials of type $\sE(o)$ where 
$|o_1 - o_2| \leq o \leq o_1 + o_2$. 

\noindent{\bf Remark.} We will see,  that these degenerated surfaces 
 are not necessary cone points of $(\sF_{\tau},\omega_{\tau})$. 
\medskip\\
Now the cylinder decomposition of $\abas \in \sF_{\tau}$ 
is not affected when deforming $\abas$ inside $\sC_i$, 
however the height of some horizontal cylinders on $\abas$ disappears 
when deforming into $\partial^{top}\sC_i$. Hence there 
are some horizontal cylinders on $\abas$ which degenerate 
to saddle connections when deforming $\abas$ into $\partial^{top}\sC_i$. 
The only thing we use is that the (width of the) cylinders on 
$\abas$ which degenerate do not depend on point 
where we enter $\partial^{top}\sC_i$ through $\sC_i$. 
\medskip\\      
This allows to simplify counting formula (\ref{rate5}) 
above using the cylinder decomposition of $\sF_{\tau}(o_1,o_2)$ 
{\em and} the $\slz$-orbit $\sO_{\alpha}$ of 
$\abas \in \sF_{\tau}(o_1,o_2)$: 
\begin{eqnarray}\label{torsioncount2}
c_{cyl}(\alpha)=
\frac{1}{|\sO_{\alpha}|}
\sum^{n_{\sF}}_{i=1}\left[ 
\sum^{n_{i}}_{k=1}\frac{|\sO_{\alpha} \cap \sC_i|}{w^2_{i,k}}+  \right. 
\left. \sum^{m_i}_{k=1}\frac{|\sO_{\alpha} \cap \partial^{top} \sC_i|}
{w^2_{i,k}}
\right].
\end{eqnarray}
{\bf Calculating quadratic constants using sequences.}
Choose a non-periodic direction $\theta \in S^1$ 
on $\abas \in \sF_{\tau}$, or equivalently on $\sF_{\tau}$,   
take the unipotent subgroup $u_{\theta}(t) \subset \slr$ 
and look at the orbit 
\[ \sO_{\theta, \alpha}(T) := 
\{u_{\theta}(t)\cdot \abas: \ t \in  (0,T) \} \]
Using the Siegel-Veech 
formula and ergodicity of the $u_{\theta}$ action on 
$\overline{\cup_{T>0} \sO_{\theta,\alpha}(T)}$, 
one can calculate all kinds of Siegel-Veech constants. 
Given the shape of the Siegel-Veech formula 
one can try to use $\sO_{\theta,\alpha}(T)$,  
to approximate Siegel-Veech constants for $\abas$ 
(however this can be done). One expects to   
be close to the real Siegel-Veech 
constant for $\abas$, when $T$ is large enough.     
However the orbit we have chosen hits 
the modular fiber $\sF_{\tau}$ only once. 
Thus one needs to consider unipotent subgroups 
which orbit returns to $\sF_{\tau}$ often, 
i.e. the unipotent subgroups generated by 
parabolic elements fixing a periodic direction 
of $\sF_{\tau}$.  Depending on the 
chosen point $\abas \in \sF_{\tau}$ we see 
a finite or infinite intersection   
 \[ \overline{\cup_{T>0} \sO_{\theta,\alpha}(T)} 
 \cap \sF_{\tau},\]
in the infinite case we see (a union of) leaves
\[ \sL = 
\overline{\cup_{T>0} \sO_{\theta,\alpha}(T) \cap 
\sF_{\tau}} 
\subset \CF_{\theta}(\sF_{\tau},\omega_{\tau}).\]

In the finite orbit case we like to ask, if 
a sequence of finite orbit surfaces 
$(S_i,\alpha_i) \in \sF_{\tau}$ 
with $\lim_{i \rightarrow \infty}|\sO_{\alpha_i}|=\infty$ 
admits the continuity property 
 \[ \lim_{i \rightarrow \infty}c_{\ast}(S_i,\alpha_i)
 =c_{\ast,gen}(S,\alpha) \]
for any type of Siegel-Veech constant $c_{\ast}$.    
Here $c_{\ast,gen}$ is the Siegel-Veech constant for 
generic surfaces in $\overline{\cup_i\sO_{\alpha_i}}$. 
Because $\sF_{\tau}$ is an arithmetic surface 
itself, we expect
\[\overline{\cup_i\sO_{\alpha_i}}=\sF_{\tau}. \]  
To avoid a general discussion,  
{\em assume} for the moment 
the orbits $\sO_{\alpha_i}$ 
become more and more equally distributed in 
$\sF_{\tau}$ with respect to Lebesque measure.   
This is known to be true for  
$\abas \in \sF_{\tau}$ with infinite $\slz$ 
orbit, in fact \cite{emm,ems} implies that 
$\sO_{\alpha}$ is equally distributed 
in $\sF_{\tau}$ with respect to Lebesque measure. 
Using either series of finite orbits 
(with equal distribution assumption), 
or infinite orbits, 
we obtain the asymptotic constant 
for {\em generic differentials} as a 
limit from formula (\ref{torsioncount2}) 
\begin{eqnarray}\label{genericcount2}
c_{cyl}(\alpha)= 
\sum^{n_{\sF}}_{i=1}
\frac{\area(\sC_i)}{\area(\sF_{\tau})}
\sum^{n_{i}}_{k=1}\frac{1}{w^2_{i,k}},
\end{eqnarray} 
because by equal distribution in the limit
\begin{equation} 
\frac{|\sO_{\alpha_i}\cap \ \sC_j|}{|\sO_{\alpha_i}|}\stackrel{i \rightarrow \infty}{\longrightarrow}
\frac{\area(\sC_j)}{\area(\sF_{\tau})}\ \text{ 
and } \
\frac{|\sO_{\alpha_i}
\cap \ \partial^{top}\sC_j|}{|\sO_{\alpha_i}|}
\stackrel{i \rightarrow \infty}{\longrightarrow}
\frac{\area(\partial^{top}\sC_j)}{\area(\sF_{\tau})}=0.
\end{equation} 
Without giving a completely general 
proof (for equal-destribution or convergence) at this place, 
we verify convergence of asymptotic constants 
in case of $d$-symmetric differentials (in most  cases).  
\vspace*{3mm}\\
{\bf Counting saddle connections.} 
The asymptotic formula for saddle connections is the 
same as the one for cylinders as long as 
$\abas \in \sF_{\tau}$ is arithmetic (or in general: Veech). 
One just needs to replace  
the width $w_i$ of cylinders in the horizontal direction by 
the length of the bounding saddle connections $s_j$ in 
formula \ref{torsioncount2}. Things become different 
if we take $\abas \in \sF_{\tau}$ and ask for the  quadratic growth 
rate of saddle connections which connect the two (different) singular points 
$z_0,z_1 \in Z(\alpha)$ of $\abas$. If a line segment  connecting 
$z_0$ with $z_1$ is horizontal, 
$\abas$ is located on a saddle connection in 
$\CF_h(\sF_{\tau})$, because we can degenerate $\abas$ 
by moving $z_1$ into $z_o$ along a saddle connection 
$s_{\sF} \in \CF_h(\sF_{\tau})$.  
Now the saddle connections on $\abas$ which are killed 
by this deformation are exactly of the length of the path 
necessary to degenerate $\abas$ along $s_{\sF} \in \CF_h(\sF_{\tau})$, 
or in other words the length of these saddle connections is the 
distance of $\abas \in s_{\sF} \subset \sF_{\tau}$ to the 
endpoint of $s_{\sF}$.  
\medskip\\
There are two possible directions (left or right) 
to move $\abas$ into a cone point of $\sF_{\tau})$ 
along the saddle connection $s_{\sF}$. 
We denote the distance of $\abas \in s_{\sF}$ to the left, or right 
endpoint of $s_{\sF}$ by $s^-_{\alpha}$, $s^+_{\alpha}$  respectively.  
On each $\abas \in s_{\sF}$ there are $m^{\pm}_s$ 
saddle connections of length 
$s^{\pm}_{\alpha}$ which disappear when moving $\abas$ 
into the endpoints on $s_{\sF}$. Using this and formula 
\ref{torsioncount2} we find the quadratic asymptotic 
constant $c_{\pm}\aby$ for the set of saddle connections 
$SC^{\pm}$ on $\aby \in \sF_{\tau}$ connecting 
two (different) cone points $z_0$ and $z_1$:  
\begin{multline}\label{torsadcount2}
c_{\pm}(\tau)=
\frac{1}{|\sO_{\tau}|}
\sum_{s \in SC_h(\sF_{\tau})}
 \sum_{\alpha \in \sO_{\tau}(s)}\left[ \frac{m^-_s}{(s^{-}_{\alpha})^2}+  
\frac{m^+_s}{(s^{+}_{\alpha})^2} 
\right]=\\ 
=\frac{2}{|\sO_{\tau}|}\sum_{s \in SC_h(\sF_{\tau})}
 \sum_{\alpha \in \sO_{\tau}(s)}\frac{m^{\pm}_s}{(s^{\pm}_{\alpha})^2}.
\end{multline}
with $\sO_{\tau}(s):=\sO_{\tau} \cap s$. The last identity follows 
from the existence of an involution $\phi \in \aff(\sF_{\tau})$ 
with $\D\phi=-\id$ acting on $SC_h(\sF_{\tau},\omega_{\tau})$.  
\medskip\\
Remarkably one can also obtain Siegel-Veech constants $c_{\pm}$ 
for generic surfaces in $\sF_{\tau}$ using  
sequences of arithmetic differentials $(S_i,\alpha_i)$. 
At hand of the modular fibers for $d$-symmetric 
differentials one sees, that the smallest pieces into which 
we can decompose the generic constant 
\[c_{\pm, gen}=\sum_{\sO_{deg}} c_{\pm}(\sO_{deg})\] 
such that  
\[ \lim_{i \rightarrow \infty}c_{\pm}(S_i,\alpha_i)=c_{\pm}(\sO_{deg})\]
are the constants associated to $\slz$ orbits of 
degenerated differentials in $\sF^c_{\tau}$. 
\medskip\\  
We show convergence of constants for 
$2$-symmetric differentials on pgs. 
\pageref{convergence} -- \pageref{convergenceend}.
Comparing our method, i.e. using sequences, with the direct calculation using the Siegel-Veech formula we observe that 
only the number of intersections  
\[\sum_{p \in \sO_q}|\sO_n \cap B_p(\varepsilon)|, \quad \frac{1}{q} 
\in Z(\omega_{\tau})\subset \sF^{sym,c}_d \]
with a small disc $B_p(\varepsilon)$ of radius $\varepsilon$, 
centered at a cone point $p$ plays a role. 
To get the right generic constant in the limit  
$n \rightarrow \infty$ $\sO_n$ must 
become more and more equally 
distributed in $\sF_{\tau}$ since only then: 
\[ \frac{1}{|\sO_n|} \sum_{p \in \sO_q}| \sO_n \cap B_p(\varepsilon)| \stackrel{n \rightarrow \infty}{\longrightarrow}
\frac{|\sO_q|}{\area(\sF_{\tau})}   \pi \varepsilon^2 
.\]
Note, that in this case it is not necessary to 
consider distances 
of a particular point or surface 
$\abas \in \sO_n \cap B_q(\varepsilon)$ 
to the cone point $q$ like in the previous formulas.  
\medskip\\
It is common, see \cite{emz} and \cite{ems}, to divide the 
generic asymptotic constant into constants by prescribing  
{\em orders of cone points} 
$o=o(z)$, $z \in Z(\omega_{\tau})$.  
\begin{equation} 
c^{\pm}_{gen}=\sum_{o=o(z),\  
z \in Z(\omega_{\tau})} c^{\pm}_{gen,o}
\end{equation} 
with 
\begin{equation} c^{\pm}_{gen,o}= 
2 \frac{\zeta (2)}{\area(\sF_{\tau})} 
 |\{z \in Z(\omega_{\tau}):o(z)=o\}| \cdot o \cdot m_o.
\end{equation}
Note that sets with prescribed orders of cone points are 
$\slz$-invariant and that our previous formulas 
represent a decomposition of constants which 
goes beyond fixing cone-point orders. 
\vspace*{2mm}\\

\noindent{\bf Other asymptotic constants.} \ 
Assume $\sF_{\tau}:=\sF_{\tau}(o_1,o_2)$ 
is a fiber of a $2$-dimensional space of 
elliptic differentials having  
two cone points of order $o_1$ and $o_2$. 
Let $\sC_1,...,\sC_n$ be the decompositon of the horizontal foliation 
$\CF_h(\sF_{\tau},\omega_{\tau})$ 
into maximal, open-cylinders. 
We might take an open cylinder $\sC_{i,(a,b)} \subset \sC_i$. 
If $t_v \in (0,h_i)$ is a transversal coordinate 
of $\sC_i$ ($h_i$ denotes the height 
of $\sC_i$), $\sC_{i,(a,b)}$ is the set of points 
in $\sC_i$ for which $t_v \in (a,b) \subset (0,h_i)$. 
Since $\sC_{i,(a,b)}$ is $u_t=\left(\begin{smallmatrix}1&t \\ 1&0 
\end{smallmatrix}\right) $ invariant we can define the following 
asymptotic constants for {\em finite orbit} 
$\abas \in \sF_{\tau}$: 
\begin{eqnarray}\label{torsioncount3}
c_{i,(a,b)}(\alpha)=
\frac{1}{|\sO_{\alpha}|}
\sum^{n_{i}}_{k=1} \frac{|\sO_{\alpha} \cap \sC_{i,(a,b)}|}{w^2_{i,k}},
\end{eqnarray}
and for {\em generic} $\abas \in \sF_{\tau}$:
\begin{eqnarray}\label{genericcount3}
c_{i,(a,b)}(\alpha)= 
\sum^{n_{\sF_{\tau}}}_{i=1}
\frac{\area(\sC_{i,(a,b)})}{\area(\sF_{\tau})}
\sum^{n_{i}}_{k=1}\frac{1}{w^2_{i,k}},
\end{eqnarray} 
The cylinders $\CC_{i,k}$ contained in the horizontal foliation 
of $\abas \in \sC_{i,(a,b)}$ have heights restricted by the condition 
$t_v \in (a,b)$. The condition also sets a restriction on the vertical 
distance of the two cone points of $\abas$. This in turn implies that 
the {\em areas} of the horizontal cylinders for $\abas \in \sC_{i,(a,b)}$ 
are restricted. Since the cylinder area is invariant under 
$\slr$ deformations the asymptotic constants 
$c_{i,(a,b)}(\alpha)$ in \ref{torsioncount3} and \ref{genericcount3}
are the asymptotic quadratic growth rates of cylinders in 
$(X,\beta) \in \sF_{\tau}$ which 
map to horizontal cylinders  
of some $\abas \in \sC_{i,(a,b)}$ for certain $A \in \slz$.  
\medskip\\   
For precise statements one can restrict to concrete examples. 
Note however that for finite orbit $\abas$ one can always choose 
intervals $(a,b)$ such that 
\[ \sO_{\alpha} \cap \sC_{i,(a,b)}= \emptyset,  \] 
in particular there are no cylinders on differentials $\aby \in \sO_{\alpha}$ 
which have area below a certain positive number.


\section{$d$-symmetric elliptic differentials}



A cone point $s$ of order $o$ on $\abx$ is called {\em degenerated}, 
if there is an open neighborhood $U$ of $s$ such that $U\backslash \{s\}$ 
is homeomorphic to the disjoint union of at least two punctured discs.
We call $s$ {\em totally degenerated} if there is a neighborhood $U$ such that 
$U\backslash \{s\}$ is homeomorphic to $o$ punctured discs. A 
differential $\abx$ is called degenerated, if it contains degenerated 
singular points. We call $\abx$ {\em totally degenerated} if {\em all} 
singular points are totally degenerated.\\

Note $U \backslash \{s\}$ is homotopy equivalent to a union of circles. 
Moreover degenerated surfaces are not described by a one form alone, 
one eventually has to provide information about identifications of 
(some) points. 
Degenerated surfaces appear naturally after collapsing 
cone points. In particular degenerated surfaces are part of the 
{\em closure of moduli space fibers} $\sF$. 
\medskip\\
{\bf $d$-symmetric torus coverings.} 
We consider elliptic differentials of the shape
\[ (X,\omega)= (\#^d_{t_h+it_v}\T^2,\#^d_{t_h+it_v}dz) \] 
with a line segment $[0,1] \cdot (t_h+it_v) \subset \C$. 
Note that $\Z/d\Z$ acts on $\#^d_{t_h+it_v}\T^2$ 
as subgroup of 
\begin{equation}
\aut(\#^d_{t_h+it_v}\T^2):=\aut(\#^d_{t_h+it_v}\T^2, \#^d_{t_h+it_v}dz),
\end{equation} 
hence the differentials $\#^d_{t_h+it_v}\T^2$ are 
$d$-symmetric. 
Recall that we have already fixed 
an ordering of the sheets $\T^2$ of  
$\pi: \#^d_{t_h+it_v}\T^2 \rightarrow \T^2$
consistent with the cut and paste along $v=t_h+it_v$ 
and the action of $\Z/d\Z$, like in 
figure \ref{symmcover} below.  

\begin{figure}[h]\label{symmcover}
\epsfxsize=10truecm
\centerline{\epsfbox{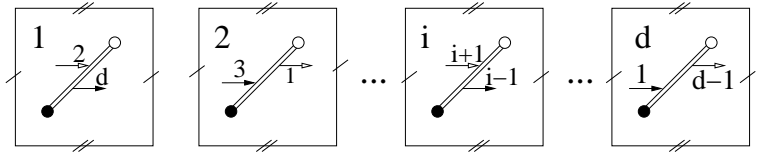}}
\caption{Canonical model of a $d$-symmetric differential}
\end{figure} 
Denote the image of the horizontal line 
$\{ \im z =s \}\subset \C$ 
on the base torus $\T^2=\C/\lat$ by $\CL_h(s)$. 
If the geodesic segment $t_h+it_v \subset \T^2$ 
is defined as the projection of the geodesic 
segment connecting $0$ and $t_h+it_v$ in $\C$, it  
intersects $\CL_h(s) \in \CF_h(\T^2)$  exactly 
\begin{itemize}
\i $[t_v]+1$ times, if $\{s\} \leq \{t_v\} $
\i $[t_v]$ times, if $\{s\} > \{t_v \}$.
\end{itemize}
Now parameterize $\CL_h(s)$ by the curve 
\[\gamma_{s}: \left\{ \begin{array}{ccc} 
[0,1]& \rightarrow &\CL_h(s) \subset \T^2\\
t & \mapsto & t + is +\lat
\end{array} \right.   
\]
Because of the cut and paste along $t_h+it_v$,  
the lift $\gamma_{s,i}$ of $\gamma_{s}$ to 
$\#^d_{t_h+it_v}\T^2$, characterized by  
$\gamma_{s,i}(0) \in \T^2_i$, is on  
the sheet $\T^2_{i+j}$ when its projection 
$\gamma_s$ has crossed the segment 
$t_h+it_v$ $j$ times. 
Given a real number $a$ we write by slight 
abuse of notation $[a]$ for the integer part 
of $a$ and $\{a\}:=a-[a]$. 
Using this we observe that 
after one loop of $\gamma_s$, 
\vspace*{1mm}\\
{\em either} \hspace*{4mm}$ \gamma_{s,i}(1)\in \T^2_{i+[t_h]+1}$, 
\ if $\{s\}  \leq \{t_h \}$\\
{\em or} \hspace*{9mm} $\gamma_{s,i}(1)\in \T^2_{i+[t_h]}$, 
\hspace*{4mm} if $\{s \} > \{t_h \}$.  
\vspace*{2mm}\\

\noindent If $\gamma_1$ and $\gamma_2$ are two chains of 
geodesic segments on an Abelian differential $\abx$ with 
\[\partial \gamma_1= \partial \gamma_2\]  
we might use cut and paste to see that 
\[ \#^d_{\gamma_1}(X,\omega)= \#^d_{\gamma_2} (X,\omega)
\] 
if the lifts of $\gamma_1$ and $\gamma_2$ to the 
universal covering $(\widetilde{X},\widetilde{\omega})$ of $\abx$ 
bound a flat disc (i.e. a disc without cone points).  
Now the line segment $t_h + i t_v \in \C$ decomposes into segments 
$[t_h]\cdot 1$, $[t_v]\cdot i$ and $\{t_h\}+i\{t_b\}$. 

\begin{lemma}\label{conepoint-structure}  
Let $\sF^{sym}_{d,\C}$ 
be the component of the fiber of elliptic 
differentials containing 
$\{\#^d_{t_h+it_v}(\T^2,dz): t_h+it_v \in \C\}$. 
Then any differential contained in  
$\sF^{sym,c}_{d,\C}-\sF^{sym}_{d,\C}$ is totally degenerated.
\end{lemma}
\begin{proof}
Note, that $\sF^{sym}_{d,\C}$ is connected, because $\C$ is. 
Take the $d$-symmetric differential 
$\#^d_{h+iv+\{t_h\}+i\{t_v\}}(\T^2,dz)$ 
associated to the line segment 
\[t_h + i t_v =[t_h]\cdot 1+[t_v]\cdot i+\{t_h\}+i\{t_v\}= 
h\cdot 1+v\cdot i+\{t_h\}+i\{t_v\}\] 
with two integers $h$ and $v$. 
Denote the $j$-th sheet of $\#^d_{h+iv+\{t_h\}+i\{t_v\}}\T^2$ 
by $\T^2_j$, $j=1,...,d$ and rewrite the line segment 
$h\cdot 1+v\cdot i+\{t_h\}+i\{t_v\}$ on each $\T^2_i$ 
as the composition of the three segments $h \cdot 1 \subset \T^2$, 
$v\cdot i \subset \T^2$ and $\{t_h\}+i\{t_v\} \subset \T^2$. 
By the discussion before the theorem we jump 
from the sheet $\T^2_j$ to the sheet $\T^2_{j\pm h (d)}$
when crossing $h \cdot 1 \subset \T^2$ {\em upwards} ($+$) or 
{\em downwards} ($-$). If we cross $v \cdot i \subset \T^2$ 
we jump from  $\T^2_j$ to the sheet $\T^2_{j\pm v (d)}$ 
depending if we cross {\em to the right} ($+$) or 
{\em to the left} ($-$). 
\medskip\\
With fixed $h,v \in \Z$, we now assume 
$\{t_h\}+i\{t_v\}=0 \in \T^2$ 
and follow a circular path $\epsilon e^{i \theta}$ 
on $\#^d_{h+iv}\T^2$ 
centered at $0 \in \T^2_j$ counter-clockwise. 
Starting with small $\theta$ we jump to the sheet 
$j-v \ (d)$ when crossing the vertical 
$i \cdot v \subset \T^2_j$ to the left. 
Next we are crossing the horizontal $h \subset \T^2_{j-v \ (d)}$   
downwards to jump onto the sheet $j-v-h \ (d)$.  
Continuing in that manner we jump back to $j-v \ (d)$ 
and then to $j \ (d)$ when completing the 
circle at $\theta = 2\pi$. This already shows 
that all degenerated $d$ symmetric covers 
$\#^d_{h+iv}(\T^2,dz)$ ($h,v \in \Z$) 
are in fact totally degenerated.   
\medskip\\
Now consider a $d$-symmetric differential defined by 
the geodesic segments $h \in \Z$, $iv \in i\Z$ 
and $\epsilon e^{i \theta}$, with $\epsilon >0$.  
Changing $\theta \in [0,2\pi)$ as above gives a 
closed loop, isometrically embedded in 
$\sF^{sym}_{d,\C}$ and shows that 
$\sF^{sym,c}_{d,\C}$ has no cone points.
\vspace*{2mm}
\end{proof}
\begin{lemma}\label{torus-structure}
Assume the moduli space $\sE^{sym}_d$ is connected, 
then the connected components of  
$\sF^{sym}_d \subset \sE^{sym}_d$ 
have no cone-points of positive order and 
$\slz \cdot\sF^{sym,c}_{d,\C}$ is isomorphic to a 
disjoint union of tori: 
\[ \slz \cdot \sF^{sym,c}_d 
\cong \slz \cdot \C/d_1\Z \oplus d_2\Z i \subset \sE^{sym,c}_d.\]  
If $\sF^{sym,c}_d$ is connected there is an 
$m \in \N$, such that: 
\[\sF^{sym,c}_d \cong \T^2_m:=\C/m\Z \oplus m\Z i.\]  
\end{lemma}
\begin{proof}
The first statement follows immediately from Lemma 
\ref{conepoint-structure}: 
all connected components of $\slz \cdot \sF^{sym,c}_d$ 
must be tori, because they are in the 
$\slz$ orbit of $\sF^{sym,c}_{d,\C} $ 
by connectedness of $\sE^{sym}_d$, hence 
affine images of the torus $\sF^{sym,c}_{d,\C} $.
In other words 
\[ \slz \cdot \sF^{sym,c}_d \cong 
\bigcup^n_{i=1}\C/ \Lambda(v_i,w_i)\]
with lattices $\Lambda(v_i,w_i)$ generated by two 
integer vectors $v_i, w_i \in \Z^2$. Because  
\[ \C/ \Lambda(v_i,w_i) \cong A_i \cdot \R^2/a\Z \oplus b\Z  \]
for some $A_i \in \slz$ all tori have the same area.  
The second statement is clear, because rotation by $90$ degrees 
is contained in $\slz$.  
\vspace*{2mm}
\end{proof}

\noindent The following Lemma  tells us that surfaces in $\sF^{sym,c}_d$ 
are actually belonging to $\sF^{sym,c}_{d,\C} $. 
\begin{lemma}\label{mod-structure}
Assume $\abx$ is $d$-symmetric and $\per(\omega)=\Z^2$.\\ 
Then there is a $t_h+ it_v \in \C$ such that 
\[ \abx \cong \#^d_{t_h+ it_v} (\T^2,dz).\]
\end{lemma}
\begin{proof}
Since $\abx$ is $d$-symmetric, we have two covering maps:
\[ \pi_d: X \rightarrow X/(\Z/d\Z)\cong \C/\Lambda \ \mbox{ and }  
\pr: X \rightarrow \C/\per(\omega)\cong \T^2.\] 
This implies there is a translation covering 
$\T^2 \rightarrow X/(\Z/d\Z)$ and by the cone point orders 
of $\abx$ we must have $X/(\Z/d\Z)\cong \T^2$.  
We assume one cone point on $\abx$, say $z_0 \in X$ 
is the preimage of $[0] \in \T^2$. If the second cone point 
is the preimage of $[a+ib] \in \T^2$ we remove the line 
segment $I_{a+ib}:=\{t[a+ib]: t \in [0,1]\} \subset \T^2$. 
We can present $X-\pi^{-1}_d(I_{a+ib})$ by $d$ squares 
\[ \CQ_{a+ib,j}:=[0,1] \times [0,i]-I_{a+ib} \subset \C \quad j=1,...,d\]
with {\em certain identifications along the boundaries} of 
$\partial \CQ_{a+ib,j}$. Further we can assume  
$\Z/d\Z$ acts by moving $\CQ_{a+ib,j}$ to $\CQ_{a+ib,j+1 \mod d}$. 
\medskip\\
Denote the four boundary line segments of $\CQ =[0,1] \times [0,i]$ 
by $\partial^{\ast} \CQ$ with $\ast=btm$ for the "bottom" 
component, $\ast=top$ for the "top", $\ast=left$ for the "left" and
$\ast=right$ for the "right" boundary component.
Since the identifications done along the boundary and along 
$I_{a+ib}$ have to match with the $\Z/d\Z$-action we must have 
\begin{itemize}
\i $\partial^{btm} \CQ_{a+ib,j} = \partial^{top} \CQ_{a+ib,j+k}$ 
\vspace*{1mm}
\i $\partial^{left} \CQ_{a+ib,j} = \partial^{right} \CQ_{a+ib,j+l}$  
\end{itemize}
for all $j=1,...,d$. 
By choice when crossing $I_{a+ib}$ 
from the left, we change from $\CQ_{a+ib,j}$ to 
$\CQ_{a+ib,j+1}$. 
Now given two numbers $k,l \in \{0,...,d\}$,  
the line segment $l+ik+a+ib \in \C$ defines 
a $d$-symmetric covering $\#^d_{(a+l)+ i(b+k)} (\T^2, dz)$ 
which has all the above properties and thus
\[ \abx \cong \#^d_{(a+l)+ i(b+k)} (\T^2,dz).\]
\end{proof}

\noindent{\bf Proof of Theorem \ref{symstructure}.}
By Lemma \ref{mod-structure} 
$\sF^{sym}_d$ is connected and thus a torus $\C/\Lambda$ 
by Lemma \ref{conepoint-structure}. The lattice $\Lambda$ 
of this torus is a sub-lattice of $\lat$ by Lemma \ref{torus-structure}, 
moreover by $\slz$-invariance of $\sF^{sym}_d$ (rotation by 
$90$ degrees!), $\Lambda = m\Z\oplus m\Z i$ for an $m \in \N$. Thus 
\[\sF^{sym}_d \cong \C/m\Z\oplus m\Z i\] 
for an $m \in \N$. 
Now for $1 > \epsilon > 0$ take the horizontal path 
\[\gamma: t \longmapsto \left[\#^d_{t+i\epsilon}\T^2 \right] 
\subset \CF_h(\sF^{sym}_d),\]  
then using Lemma \ref{mod-structure} 
it is easy to see that $\gamma(0)=\gamma(d)$.  
One can argue as follows: split the line segment $d+i\epsilon$ 
into two segments $d$ and $i\epsilon$ 
and consider the 
surface $\#^d_{d+i\epsilon}\T^2$. The  
monodromy of $\#^d_{d+i\epsilon}\T^2$ 
agrees with the monodromy of 
$\#^d_{i\epsilon}\T^2=\#^d_{0+i\epsilon}\T^2$ 
and thus  
$\#^d_{d+i\epsilon}\T^2 \cong \#^d_{0+i\epsilon}\T^2$, 
hence $m=d$.    
\qed 
\medskip\\
{\bf Twist coordinates on $\T^2_d$.} 
For $1> t_v >0$ consider the 
$d$-symmetric differential $\#^d_{t_h+it_v}\T^2$. 
Then the horizontal foliation of 
$\#^d_{t_h+it_v}\T^2$ contains $d+1$ 
cylinders: $d$ of width $1$ attached to the top of one 
cylinder of width $d$.  
We call $t_h$ the {\em horizontal twist} and $t_v$ the {\em vertical twist} 
of $\#^d_{t_h+it_v}\T^2$. We just found, that 
parameterization of $d$-symmetric differentials by horizontal- and 
vertical twists gives an isometric universal covering
\begin{equation}
\label{hvparameter}
\begin{CD}
\C  \longrightarrow  \T^2_d\\
t_h+it_v \longmapsto  \left[\#^d_{t_h+it_v}\T^2 \right]   \\
\end{CD}
\end{equation}
\vspace*{2mm}


\noindent {\bf Global description of ${\mathbf \sE^{sym}_d}$.} To describe 
the space $\sE^{sym}_d$ of all $d$-symmetric 
differentials (with fixed area $d$), we note that there is 
a fibration  
\begin{equation}
\sF^{sym}_d \longrightarrow \sE^{sym}_d \longrightarrow \slr/\slz.
\end{equation}
Since $\sF^{sym}_d \cong \T^2_d$, we can describe 
$\sE^{sym}_d$ as 
\begin{equation}
\sE^{sym}_d \cong \slr \ltimes \C_{\lat} /\slz \ltimes d(\lat).
\end{equation}
with $\C_{\lat}:= \C - \lat$.
Note: for all $d$-symmetric differentials $\#^d_{t_h+it_v}\T^2$ there 
is an involution $\phi \in \aff(\#^d_{t_h+it_v}\T^2)$, i.e. $\phi$ is 
affine linear with $\D\phi =-\id$.  $\phi$ exchanges the 
two cone-points of order $d$ on $\#^d_{t_h+it_v}\T^2$. 
Consequently classifying $d$-symmetric differentials without 
{\em named} cone-points gives a sphere (lattice points removed) 
$\proj^1_{d,\lat} := \sF^{sym}_d/(-\id)=\T^2_{d,\lat}/(-\id)$, 
which admits only a {\em quadratic differential}. 
The $\slr$ action on $\sF^{sym}_d$ descends to 
an $\pslr$ action on 
$\proj^1_{d,\lat}$ and the moduli space $\sE^{sym}_{d,\pm}$ 
has the following fiber bundle structure:  
\begin{equation}
\proj^1_{d,\lat} \longrightarrow  \sE^{sym}_{d,\pm} 
\longrightarrow \pslr/\pslz
\end{equation}
Forgetting the elliptic differential structure of 
$(\#^d_{t_h+it_v}\T^2,\#^d_{t_h+it_v}dz)$ gives the 
moduli space $\sM_d^{sym}$ $d$-symmetric torus covers of 
genus $d$ with two branch points of order $d-1$. 
Algebraically that is dividing out the circle  
bundle defined by the action of 
$\so$ on $\sE^{sym}_d$    
\begin{equation}
\proj^1_{d,\lat} \longrightarrow \sM_d^{sym} \longrightarrow \H/\pslv{\Z\
}.
\end{equation}


\section{Geometric properties of $d$-symmetric  differentials} 


Now we compose our knowledge about the modular fiber $\T^2_d$ 
of $d$ symmetric torus covers over the base lattice $\lat$ 
and our counting formula to find the quadratic growth rate 
function on $\T^2_d$. 
\medskip\\

 
\begin{prop}\label{horfol}
Take the differential $\#^d_{t_h+it_v}\T^2$  
and assume $0<t_h<1$. Then the horizontal 
foliation of $\#^d_{t_h+it_v}\T^2$
contains  
\begin{itemize}  
\i $([t_v],d)$  cylinders of width $\frac{d}{([t_v],d)}$ and 
\i $([t_v]+1,d)$ cylinders of width $\frac{d}{([t_v]+1,d)}$, 
if $t_v \notin \Z$ \vspace*{2mm}
\i $([t_v],d)$  cylinders of width $\frac{d}{([t_v],d)}$, if $t_v \in \Z$.
\end{itemize}
We use $(0,d)=d$ here.
\end{prop}
\begin{proof}
By the discussion before the Proposition 
the point $\gamma_{i,s}(1)$ 
is in  $\T^2_{i+[t_v]+1}$, if $s \mod 1 \leq t_h \mod 1$ and 
in $\T^2_{i+[t_v]}$, if $s \mod 1 > t_h \mod 1$. 
By $\Z/d\Z$ symmetry continuing $\gamma_{i,s}$ 
we are back on the leaf $i$ as soon as 
$n([t_v]+1) \equiv 0 \mod d$ in the 
first case and $n \cdot [t_v] \equiv 0 \mod d$ in the second.
That means the cylinder $\sC$ containing $\gamma_{i,s}$ has 
width $w=\frac{d}{([t_v]+1,d)}$ if $s \mod 1 \leq t_v \mod 1$ and 
width $w=\frac{d}{([t_v],d)}$ if $s \mod 1 < t_v \mod 1$. 
In case $([t_v]+1,d)>1$ or $([t_v],d)>1$ there are 
sheets $\T^2_k$ which $\sC$ does not touch. 
Using the $\Z/d\Z$ action we find all other horizontal cylinders 
containing the preimages of the curve $\gamma_s$.  
Since the total length of all the preimages of $\gamma_s$ is 
$d$, there are exactly $([t_v]+1,d)$ 
cylinders $\sC$ above $\gamma_s$ in the first case 
and $([t_v],d)$ in the second. 
\end{proof}
The horizontal foliation of $\sF^{sym}_d$ decomposes into $d$ 
open cylinders $\sC_0,...,\sC_{d-1}$, each of width $d$ and height $1$.
The surface $\#^d_{t_h+it_v}\T^2$ belongs to $\sC_i$, 
if and only if $[t_v]=i$ and $t_v \notin \Z$. 
If $t_v =i \in \Z$, 
\[\#^d_{t_h+it_v}\T^2 \in \partial^{top} \sC_{i-1}
=\partial^{btm} \sC_{i}.\]
\begin{theorem}\label{decomp}
The horizontal foliation of $\abx \in \sF^{sym}_d$ 
contains 
\begin{itemize}  
\i $(i,d)$  cylinders of width $\frac{d}{(i,d)}$ and \vspace*{1mm}
\i $(i+1,d)$ cylinders of width $\frac{d}{(i+1,d)}$, 
\hspace*{2mm}if $\abx \in \sC_i$ \vspace*{1mm}
\i $(i,d)$  cylinders of width $\frac{d}{(i,d)}$, \hspace*{2mm} if 
$\abx \in \partial^{btm} \sC_i$.
\end{itemize}
\end{theorem}
\begin{proof}
The number and width of 
closed cylinders in the horizontal foliation of 
$X=\#^d_{t_h+it_v}\T^2$ are the same for all 
$X \in \sC_i \subset \CF_h(\sF^{sym}_d)$. 
The statement now follows from Proposition \ref{horfol}. 
\end{proof}
The vertical foliation $\CF_v(\sF^{sym}_d)$ 
admits a decomposition into open cylinders $\sC_{v,0},...,\sC_{v,d-1}$ 
(of width $d$ and height $1$) and boundary components 
$\partial^{btm}\sC_{v,0},...$ $...,\partial^{btm}\sC_{v,d-1}$ 
as well and we have:  
\begin{kor}\label{vdecomp}
The vertical foliation of $\abx \in \sF^{sym}_d$ contains
\begin{itemize}
\i $(i,d)$  cylinders of width $\frac{d}{(i,d)}$ and \vspace*{1mm} 
\i $(i+1,d)$ cylinders of width $\frac{d}{(i+1,d)}$, 
if $\abx \in \sC_i$ \vspace*{1mm}
\i $(i,d)$  cylinders of width $\frac{d}{(i,d)}$, if 
$\abx \in \partial^{btm} \sC_i$.
\end{itemize}   
\end{kor}
\begin{proof}To prove this statement, we rotate the whole space 
$\T^2_d = \sF^{sym}_d$  by
\[r_{\pi/2}:=\left[\begin{smallmatrix}0&-1 \\ 1&0  
\end{smallmatrix}\right] \in \slz.
\] 
That maps the surface $\#^d_{t_h+it_v}\T^2 \in \T^2_d$ 
to 
\[r_{\pi/2} \cdot \#^d_{t_h+it_v}\T^2 = 
\#^d_{t_v-it_h}\T^2 
\in \T^2_d\] 
and isometrically maps the horizontal foliation of 
$\#^d_{t_h+it_v}\T^2$ to the vertical foliation of 
$\#^d_{t_v-it_h}\T^2$. On the level of  $\sF^{sym}_d$ 
$r_{\pi/2}$ maps the horizontal cylinder $\sC_i$ to $\sC_{v,d-i}$. 
Now notice that differentials in $\sC_i$ and $\sC_{d-i}$ 
just differ by the (hyperelliptic) involution $-\id$.
That proves the claim.
\end{proof}
Note that for prime $d>2$ there are only $2$ possible numbers 
of {\em horizontal cylinders} on differentials in 
the $\sC_i$, namely $2$ and $d+1$.
\medskip\\    
{\bf Orbits of torsion points.} If 
the {\em prime-factor decomposition} of $n \in \Z$ is
$n=p^{l_1}_1 \cdot p^{l_2}_2 \dots p^{l_r}_r$, 
then the number of positive {\em divisors} $D(n)$ 
of $n$ is 
\[D(n)=(l_1+1)\dots(l_r+1)\] 
and we have 
\begin{prop}\label{splitting}
The set $\T^2_d[m]$ on 
$\T^2_d \cong \sF^{sym, c}_d$ contains precisely $D(m)$ 
$\slz$-orbits. In particular the set of degenerated differentials, 
represented by $\T^2_d[d] \subset \sF^{sym,c}_d$, is the disjoint union of  $D(d)$ $\slz$-orbits.  
\end{prop}
We call the periodic foliation 
$\CF_{\theta}(\#^d_v \T^2)$ 
of a $d$-symmetric differential $\#^d_v \T^2$ {\em simple}
if every cylinder $\sC_j$ in  $\CF_{\theta}(\#^d_v \T^2)$ 
is in the $\Z/d\Z$ orbit of any other cylinder $\sC_i$. 
\begin{theorem}
For $v \in \T^2_d \cong \sF^{sym}_d$ the following 
statements are equivalent: 
\begin{itemize}
\inum{1.} The $d$-symmetric differential $\#^d_v \T^2$ contains 
a simple direction with precisely $k|d$ closed cylinders.  
\inum{2.} There is an $A \in \slz$, and an $m \in \Z/d\Z$, 
such that 
\[ [A \cdot v]\in \partial^{top}\sC_{m} \mbox{ with } k=(m,d).\]
\end{itemize}
Moreover, if $v \in \T^2_d$ is a torsion point, we can 
add the statement: 
\begin{itemize}
\inum{3.} There is an $n \in \N$ with $\frac{d}{(d,m)}|n$, 
 and $v \in \T^2_d(n)$. 
\end{itemize}
\end{theorem}
\begin{proof}
By Theorem \ref{decomp} the horizontal foliation of 
a $d$-symmetric surface $\#^d_v \T^2$ contains precisely 
$k|d$ cylinders 
if $v \in \partial^{top}\sC_{m}$ with $k=(m,d)$. 
Now if a simple, periodic foliation of $\#^d_v \T^2$ 
contains precisely $k|d$ cylinders there is an $A \in \slz$ 
making this direction horizontal and thus we must have 
$A\cdot v \in \partial^{top}\sC_{m}$ with $k=(m,d)$.  

Finally we need to find $A \in \slz$ such that
\[ \left[\begin{smallmatrix} \ast \\ m/d  \end{smallmatrix}\right]= 
A\cdot \left[\begin{smallmatrix} 1/n \\ 0 \end{smallmatrix}\right] \quad 
\mbox{on } \T^2. 
\]
Note that $k=(m,d)|d$. Since $\sO_n=\slz \cdot [1/n,0]=\T^2(n)$, 
we have the condition $d/k=d/(m,d)|n$.
\end{proof}
\noindent{\em Remark:} 
Suppose $\sF$ is a $2$-dimensional 
modular fiber of $\T^2$ coverings and the natural projection 
$\pi_{\ast}: \sF \rightarrow \T^2$ factors over 
$\T^2_d$, then Proposition \ref{splitting} implies that 
the preimage 
\[\sF[d]:=\pi^{-1}_{\ast}(\T^2_d[d])\subset \sF \] 
of the $d$-torsion points on $\T^2_d$ splits into 
at least $D(d)$ different $\slz$ orbits. If $\sF$ is a 
space of differentials with two cone points,  
the origamis, i.e. surfaces which can be tiled 
by unit squares, in $\sF$ are a subset of $\sF[d]$. 
\smallskip\\
We provide examples for this factorization of $\pi_{\ast}$ 
in \cite{s2,s3}. However for $\sF_d(2) \subset \sF_d(1,1)$ 
the splitting of orbits is already known, see \cite{mcm4,hl}. 
\bigskip\\
{\bf Differentials on integer coordinates.} 
We like to describe the degenerated surfaces having  
integer twist coordinates, $(j,k) \in \Z^2/d\Z^2 \subset 
\sF^{sym,c}_d$. For the first we do not care about, 
whether some of these degenerated surfaces are isomorphic, we 
just study their shape and their $\slz$-orbits in 
$\sF^{sym,c}_d$.      

\begin{theorem}\label{degsurf} For $(j,k) \in \Z^2/d\Z^2$ the degenerated 
surface $\#^d_{j+ik}\T^2$ is a union of 
$(j,k,d) :=\gcd(j,k,d)$ tori, with all integer lattice 
points identified. Each torus-component contained 
in $\#^d_{j+ik}\T^2$ has area $d/(j,k,d)$ and contains 
$(k,d)/(j,k,d)$ horizontal cylinders of width 
$d/(k,d)$ and $(j,d)/(j,k,d)$ vertical cylinders of width 
$d/(j,d)$.
\end{theorem}
\begin{proof}
Look at the surfaces with twist coordinate $j+i0 \in \Z^2/d\Z^2$. 
The surface $\#^d_{j}\T^2$ parameterized by 
$j \in \Z/d\Z$ has $d$ cylinders 
of width $1$ in the horizontal direction and $(j,d)$ 
vertical cylinders of width $d/(j,d)$. Since $\#^d_{j}\T^2$ is a one point 
union of tori, it consists of $(j,d)$ tori each of width $1$ 
and height $d/(j,d)$. Since every differential contained in 
$\Z^2/d\Z^2$ is in the $\slz$ orbit of one of the surfaces $\#^d_{j}\T^2$,  
$j \in \Z/d\Z$ and we know all the degenerated surfaces  
$\#^d_{j+ik}\T^2$ are on the $\slz$ orbit of 
$\#^d_{a}\T^2$ if and only if $(j,k,d)=(a,d)$, all surfaces $\#^d_{j+ik}\T^2$ 
with fixed $(j,k,d)$ are one point unions of $(j,k,d)$ tori, 
each of area $d/(j,k,d)$. 
\medskip\\
The second statement now follows from the fact that $\#^d_{j+ik}\T^2$ 
admits $(k,d)$ horizontal cylinders of width $d/(k,d)$ and 
$(j,d)$ vertical cylinders of width $d/(j,d)$ shown in Theorem 
\ref{decomp} and Corollary \ref{vdecomp}.\vspace*{2mm}
\end{proof}
\noindent{\bf Isomorphic surfaces on the integer lattice.}
Assume $\#^d_{j+ik}\T^2 \in \Z^2/d\Z^2 \in \sF^{sym,c}_d$, 
then 
\[\#^d_{j+ik}\T^2 \cong \#^d_{-j-ik}\T^2 = r_{\pi}\cdot\#^d_{j+ik}\T^2\]
are isomorphic. Mainly because surfaces represented by the integer 
lattice $\Z^2/d\Z^2$ admit only one cone point and thus the 
property of distinguishing cone points, which allow to distinguish 
$\#^d_{j+ik}\T^2$ and $r_{\pi}\cdot \#^d_{j+ik}\T^2$ disappears. 
\medskip\\
{\bf Splitting of orbits -- moduli space recognition of   invariant.}  
Recently Hubert and Lelievre \cite{hl} followed by a generalized 
approach of McMullen \cite{mcm4} have described invariants to distinguish 
$\slz$-orbits of genus $2$ elliptic differentials with a single zero 
(of order $2$). We suspect that one can observe a similar phenomen
in a more general setup. In 
particular for differentials in genus $2$ with $2$ zeros of order $1$ 
\cite{s3,s4}.


\section{Asymptotic constants --- cylinders, generic case} 

{\bf Siegel Veech constants for generic differentials.} 
We write down the Siegel-Veech constants as function of 
the point in the fiber of the moduli space.  
As before we restrict our considerations 
to the fiber $\sF^{sym}_d$.
\medskip\\
{\bf Proof of Theorem \ref{siegelveechdsymm}.}
Because $\abx$ in not arithmetic, $\abx \in \sF^{sym}_d=\T^2_d$ 
is not a torsion point and we apply Theorem 
\ref{siegelveechdsymm}, Formula \ref{genericcount} to obtain the result.

Now the horizontal foliation of $\sF^{sym}_d$ decomposes into $d$ 
open cylinders $\sC_1,...,\sC_d$, each of width $d$ and height $1$. 
We can label the cylinders by the (integer) vertical twist $i=[t_v]$. 
By Theorem \ref{decomp} we have $(i+1,d)$ cylinders of width $d/(i+1,d)$ 
and $(i,d)$ cylinders of length $d/(i,d)$. Cases when the vertical 
twist $t_v$ is integer play no role for the generic constant. 
Now $i=[t_v]$ is a natural number ranging from $1$ to $d$. 
Since each strip with integer vertical twist $i$ has area $d$ 
and $\area(\sF^{sym}_d)=d^2$, we find the generic 
Siegel-Veech constant for cylinders of closed geodesics: 
\begin{multline}
c_{gen}(S)=\sum^{d}_{i=1}\frac{\area(\sC_{i})}{\area(\sF^{sym}_d)} 
\left[\frac {(i,d)^3}{d^2}+\frac {(i+1,d)^3}{d^2}\right]
=\frac{2}{d^3}\sum^{d}_{i=1} (i,d)^3.
\end{multline} 
Now for given $p|d$ the number of $\overline{i} \in \Z/d\Z$ with $p=(i,d)$ 
is $\varphi(d/p)$ and thus:   
\begin{equation}
c_{gen}(S)=\frac{2}{d^3}\sum^{d}_{i=1} (i,d)^3= 
2\sum_{p|d}\frac{\varphi(p)}{p^3}.
\end{equation}
\hfill \qed\\
{\bf Limit genus to $\mathbf{\infty}$}. 
For genus to infinity we make the following 
simple limit consideration:    
\[\frac{c_{gen}}{d}=\frac{2}{d}\sum_{p|d}\frac{\varphi(p)}{p^3} \leq 
\frac{2\zeta(2)}{d} \stackrel{d \rightarrow \infty}{\longrightarrow}0.\]  
Of course $d$-symmetric differentials contain no unfolded rational 
billiard tables and it is hard to start speculations how their 
growth rates behave for genus to infinity based on $d$-symmetric 
differentials. 
A second remark: the infinite genus surface 
$\#^{\infty}_{v}\T^2$ has {\em infinitely many} 
cylinders of width $1$ in the vertical and horizontal direction. 
Thus one cannot even start a reasonable direct counting of 
closed cylinders. 


\section{Asymptotic constants --- cylinders, finite orbit case}

What is left is to calculate the Siegel-Veech constants 
for $d$-symmetric differentials represented by torsion points of 
order $n$ in $\T^2_d$. For fixed $n\in \N$ 
take a natural number $1 \leq a \leq n$ and define
\[ \frac{\pi}{\zeta(2)}c_{d,n}(a):= 
\lim_{T \rightarrow \infty }\frac{N(V_n(a), T)}{T^2}\]
to be the quadratic growth rate of the distribution
\begin{equation}
\begin{split}
V_n(a):=\{A \cdot \hol(c): c \in  
Cyl(\CF_h(\#^d_{v}\T^2)) \mbox{ with }\\ v \in \sL_{ad/n} 
\cap \T^2_d(n), \ 
A \in \slz\}.
\end{split}
\end{equation} 
\begin{lemma}\label{torsionformula}
For $1 \leq a \leq n$ and $t_v(a):=[ad/n]\in \Z$ we have 
\begin{equation}
c_{d,n}(a)= \frac{n}{\varphi(n)\psi(n)} 
 \frac{\varphi((a,n))}{(a,n)}\left(\frac{(t_v(a),d)^3}{d^2}+
\frac{(t_v(a)+1,d)^3}{d^2} \right).
\end{equation}
\end{lemma}
\begin{proof}
The number of torsion points $\T^2_d(n)$ of order $n$ on  
the horizontal leaf $\sL_{ad/n} \subset \T^2_d$ 
going through the point $[i ad/n] \in \T^2_d$ 
equals 
\begin{multline*}
|\{b \in \Z/n\Z: (b,a,n)=1 \}|=| \{b \in \Z/n\Z: (b,(a,n))=1 \}|=\\
=\begin{cases}
|\psi^{-1}((\Z/(a,n)\Z)^{\ast})|, & \text{ if } (a,n)\geq 2 \vspace*{1mm}\\ 
n, & \text{ if } (a,n)=1
\end{cases}
\end{multline*} 
where $\psi:\Z/n\Z \rightarrow \Z/(a,n)\Z$ is the natural 
homomorphism. Thus 
\begin{equation} \label{orbitnumber} |\sL_{ad/n} \cap \T^2_d(n)|=
n\frac{\varphi((a,n))}{(a,n)}.\end{equation} 
By Theorem (\ref{decomp}) the horizontal foliation of a differential 
$\#^d_w\T^2 \in \sL_{ad/n} \cap \T^2_d(n)$ 
\begin{itemize}
\i always contains $(t_v(a),d)$ cylinders of 
width $d/(t_v(a),d)$ \vspace*{1mm}
\i and contains $(t_v(a)+1,d)$ cylinders of width $d/(t_v(a)+1,d)$ 
if $\ t_v(a)\notin \Z$. 
\vspace*{1mm}
\end{itemize}
With $|\T^2_d(n)|=\varphi(n)\psi(n)$ we find the quadratic growth rates 
above. 
\end{proof}
{\bf Remark.} The first part of the argument in the 
proof of the Lemma tells us that the Teichm\"uller disk 
through the torsion points $\T^2(n)$ of order $n$  
has 
\begin{equation}
cu(n)=\frac{1}{2}\sum^n_{a=1} \varphi((a,n))=\frac{1}{2}
\sum_{l|n}\varphi \left( \frac{n}{l}\right)\varphi(l) 
\end{equation}
cusps if $n \geq 3$, and $2$ cusps if $n=2$.
\medskip\\
Now the Siegel-Veech constant $c_{d,n}= \sum^n_{a=1}c_{d,n}(a)$ 
for periodic cylinders on $d$-symmetric differentials 
$\#^d_w\T^2$ with $[w] \in \T^2_d(n)$ is
\begin{equation}
\begin{split}
c_{d,n}= \frac{n}{\varphi(n)\psi(n)}
\left(
\sum^n_{a=1} 
\frac{\varphi((a,n))}{(a,n)}\frac{(t_v(a),d)^3}{d^2}\right. + \\ 
+ \left. \sum_{ ad/n \notin \Z}
\frac{\varphi((a,n))}{(a,n)}
\frac{(t_v(a)+1,d)^3}{d^2}
\right).
\end{split}\vspace*{2mm}
\end{equation} 
To simplify this expression we {\bf restrict to torsion points of 
order ${\mathbf n}$ with $\mathbf{(n,d)=1}$}. 
Then all numbers $a$ such that $\overline{ad/n} \in \Z/d\Z$ 
are multiples of $n$, consequently 
\begin{equation}\label{baseline}
|\T^2_d(n) \cap \partial^{btm}\sC_i| =
\begin{cases} 
\varphi(n)& \text{if } i=0, \\
0& \text{if } i \neq 0.
\end{cases}
\end{equation} 
and: 
\begin{equation}
c_{d}(n)=\frac{n}{\varphi(n)\psi(n)}
\left( d\frac{\varphi(n)}{n}+
2\sum_{a\not\equiv 0 \mod n} 
 \frac{\varphi((a,n))}{(a,n)}\frac{(t_v(a),d)^3}{d^2}
\right).\vspace*{2mm}
\end{equation}
We always assume $a \in \{-n+1,-n+2,...-1,0,1,...,n-1\}$,  
representing a class $\bar{a} \in \Z/n\Z$. 
If furthermore $\mathbf{d}$ {\bf is prime} the numbers of isotopy 
classes in a given direction is 
either $1,2,d$ or $d+1$ and the quadratic growth rates are: 
\smallskip\\
\begin{center}\label{table}
\begin{tabular}{|c | c |}
\hline 
$\#$ of cylinders &  $c_d(n)=$ \\ 
\hline 
& \\
1 & $\frac{n}{\varphi(n)\psi(n)}  
\frac{1}{d^2}\sum_{d a/n \in \{1,...,d-1\}}\frac{\varphi((a,n))}{(a,n)}$ \\
& \\
2 & $\frac{n}{\varphi(n)\psi(n)}
\frac{2}{d^2}\sum_{d-1>\frac{a}{n}d>1}\frac{\varphi((a,n))}{(a,n)}$ \\
& \\
d & $\frac{n}{\varphi(n)\psi(n)}\frac{d\varphi(n)}{n}=
\frac{d}{\psi(n)}$\\
& \\
d+1 & $\frac{n}{\varphi(n)\psi(n)}
2\left(d+ \frac{1}{d^2}\right)\sum_{0< |a| < n/d}
\frac{\varphi((a,n))}{(a,n)}$ \\
& \\
\hline
\end{tabular}\vspace*{3mm}
\end{center}
In particular there are two possibilities for prime $d$:\vspace*{1mm}\\
{\bf--- either} $(n,d)=1$, then there are {\em no} directions 
with only {\em one isotopy class} on surfaces located on the $\slz$ 
orbit of $[1/n]$ and we find  
\begin{equation}\label{primeone}
\begin{split}c_d(n)=
\frac{n}{\varphi(n)\psi(n)}   \left( 
\frac{2}{d^2}\sum_{1< \frac{a}{n}d < d-1}\frac{\varphi((a,n))}{(a,n)}+
\frac{d\varphi(n)}{n}\right. + \\ + 
\left. 2\left(d+ \frac{1}{d^2}\right)\sum_{0< \frac{|a|}{n}d < 1} 
\frac{\varphi((a,n))}{(a,n)}
  \right),
\end{split}\vspace*{1mm}
\end{equation}
{\bf--- or} $d$ divides $n$, then the asymptotic constant is the sum 
of all four terms in the table above, in particular there are 
{\em always} directions with only {\em one isotopy class} on 
surfaces in the orbit of $[1/n]$.\\ 

\noindent {\bf Results for $d=2$.} For small $d$ the above formulas 
simplify. In particular for $d=2$ and $n$ {\em odd},  
formula (\ref{primeone}) reads 
\begin{equation}
\begin{split}
c_2(n)=
\frac{1}{\varphi(n)\psi(n)}   \left(2\varphi(n)+ 
\frac{9n}{4}\sum_{ a \not\equiv 0 \mod n} \frac{\varphi((a,n))}{(a,n)}
\right)= \\ 
 =   
\frac{2}{\psi(n)} + \frac{9}{4}\left(1-\frac{1}{\psi(n)}\right)=
\frac{9}{4}- \frac{1}{4\psi(n)}.
\end{split}
\end{equation}
The only thing we have used here is the obvious identity 
derived from equation (\ref{orbitnumber}):
\[ \varphi(n)\psi(n)=n\sum^n_{a=1} 
\frac{\varphi((a,n))}{(a,n)}=\varphi(n)+n\sum^{n-1}_{a=1} 
\frac{\varphi((a,n))}{(a,n)}. \]
Note, if $d=2$ the $d$-cylinder directions 
are in fact $2$-cylinder directions. Because $1 \leq a \leq n$  
for a system of representants $a$ modulo $n$ there is no $|a|>n/2$.  
In more geometrical terms: the modular fiber has two horizontal cylinders, 
and the horizontal foliation of surfaces in both of 
these cylinders have $3=2+1$ isotopy-classes of closed geodesics.  
\medskip\\
For {\em even} $n$ the twist $a=n/d =n/2$ is integer 
and is on the horizontal leaf with vertical twist one. 
This horizontal leaf (see table above) reflects all
directions with exactly one closed cylinder,    
thus for {\em even} $n$ we obtain: 
 
\begin{equation}
\begin{split}
c_2(n)&=
\frac{1}{\varphi(n)\psi(n)}   \left(2\varphi(n)+ 
 \frac{1}{4}\varphi(\frac{n}{2}) +
\frac{9n}{2}\sum_{2a/n \neq 1,2} \frac{\varphi((a,n))}{(a,n)}
\right) = \\
&=\begin{cases}   
\frac{9}{4\psi(n)} + \frac{9}{4}\left(1-\frac{2}{\psi(n)}\right)\ =
\frac{9}{4}- \frac{9}{4\psi(n)}&
 \text{ if } 4  \nmid  n \vspace*{1mm}\\
\frac{17}{8\psi(n)} + \frac{9}{4}\left(1-\frac{3}{2\psi(n)}\right)=
\frac{9}{4}- \frac{5}{4\psi(n)}
& \text{ if } 4\  | \ n.
\end{cases}
\end{split}
\end{equation}
We have used
\begin{equation}
\varphi\left(\frac{n}{2}\right)=\left \{ 
\begin{array}{cc}
\ \varphi(n), &\mbox{ if } 4  \nmid  n \\
\frac{1}{2}\varphi(n), &\mbox{ if }   4 \ | \ n.
\end{array}
\right. 
\end{equation}
\vspace*{3mm}


\section{Asymptotic constants --- saddle connections} 


Since $d$-symmetric differentials have automorphism group 
$\Z/d\Z$ and are coverings of degree $d$, 
they are are {\em balanced}. 
Take $\#^d_{t_h+it_v}\T^2$ $d$-symmetric and let 
$\pi: \#^d_{t_h+it_v}\T^2 \rightarrow (\T^2,[0],[t_h+it_v])$ 
be the canonical projection, 
then any saddle connection $s \in SC(\T^2,[0],[t_h+it_v])$ has 
$d$ preimages 
\[\pi^{-1}(s)=\{s_1,...,s_d\}=\{\aut(\#^d_{t_h+it_v}\T^2) \cdot s_1\} 
\subset SC(\#^d_{t_h+it_v}\T^2).\] 
Thus to find the quadratic growth rate of all saddle connections on 
$\#^d_{t_h+it_v}\T^2$ is easy, if we know the one for the 
two marked torus $(\T^2,[0],[t_h+it_v])$. We did this in \cite{s1}, 
however it is much easier to use the method developed in this 
paper. The main thing to do is to find the asymptotic quadratic growth 
constants for saddle connections on $\T^2$ connecting $[0]$ and $[t_h+it_v]$. 
Assuming $[t_h+it_v] \in \T^2(n)=\slz \cdot [1/n]$ formula \ref{torsadcount} 
says mainly we need to measure distance of 
the points in $\slz \cdot [1/n] \cap \sL_0$ to the origin to the 
left and right and sum up the squares of their reciprocals. 
This gives 
\begin{equation}
c_{\pm}(n)=\frac{2n^2}{\varphi(n)\psi(n)}\sum_{(i,n)=1}\frac{1}{i^2}
\end{equation}
with limit $n \rightarrow \infty$ the Siegel-Veech constant 
for saddle connections connecting the two markings on 
the generic $2$-marked torus: 
\begin{equation}
c_{\pm}=\frac{\pi^2}{3}.
\end{equation}
We obtain the obvious 
\begin{kor}
The asymptotic quadratic growth rate $c^d_{\pm}(n)$ of saddle connections 
$SC_{\pm}(\#^d_{t_h+it_v}\T^2)$ on  
$\#^d_{t_h+it_v}\T^2 \in \sO_n \in \sF^{sym}_d$  connecting the two cone 
points of $\#^d_{t_h+it_v}\T^2$ is 
\begin{equation}
c^d_{\pm}(n)=d\frac{2n^2}{\varphi(n)\psi(n)}\sum_{(i,n)=1}\frac{1}{i^2}.
\end{equation}
For infinite orbit surfaces (generic) we obtain: 
\begin{equation}
c^d_{\pm}=d\frac{\pi^2}{3}.
\end{equation}
\end{kor} 
This is clearly not surprising, but one can ask for certain 
subsets of this set of saddle connections, 
which are connected to the number of connected components  a degenerate surface in 
$\Z^2/d\Z^2 \subset \sF^{sym,c}_d$ contains. 
By $\slz$ invariance of the asymptotic constants 
we need to consider the  $\slz$ orbits on 
$\Z^2/d\Z^2$. As we saw earlier for $m|d$ 
each orbit $\sO_{d/m}$ is 
generated by a particular torus 
$\#_{m} \T^2 \in \Z^2/d\Z^2$  
and hence these orbits are in 
one to one correspondence  
to the divisors $m | d$. It follows from 
Theorem \ref{degsurf} that 
all degenerated surfaces associated to $m|d$ 
have $m$ connected components after 
removing their "cone" point. Moreover all 
these connected components are tiled by 
$d/m$ squares and contain an equal amount  
of horizontal saddle connections. 
Now we use 
this observation to define subsets of 
all saddle connections connecting the $2$ cone points of any $\abas \in \sF^{sym}_d$. 

To begin with take the (horizontal) spine  
$\sS\!\sP_h \subset \CF_h(\sF^{sym}_d,\omega_d)$ 
and look at the saddle connections $\sS\!\sP_{h,m} 
\subset \sS\!\sP_h$ emanating or terminating 
in one of the points in $\sO_m$. 
A surface $\abas$ parameterized by 
$\mathscr{S\!P}_{h,m} \subset \sF^{sym}_d$ degenerates into $\sO_m$ along $\sS\!\sP_{h,m}$. 
During such a deformation $d$ horizontal saddle connections of $\abas$ connecting the 
two cone points will degenerate.   
Since by Theorem (\ref{degsurf}) the degenerated 
surface has $m$ 
connected components, this component  
can be "cut out" of $\abas$, by taking $d/m$ 
of the $d$ saddle connections away. 
As a consequence we can 
organize these $d$ saddle connections 
into  $m$ chains of $d/m$ saddle connections 
which altogether define a boundary in homology. 
Given $d$ and $m|d$, we call such a 
configuration of saddle connections 
an $m$-{\em homologous family}.
\medskip\\ 
Note that this notion is well-defined for 
non-arithmetic (= generic) surfaces in 
$\sF^{sym}_d$.    
In particular the Siegel-Veech constant 
for $m$-homologous families (one chain)  
on the generic surface in $\sF^{sym}_d$ is    
\begin{equation}
c^d_{\pm}(m)=\frac{\pi^2}{3}\frac{\varphi(d/m)\psi(d/m)}{d \cdot m}.
\end{equation}
\medskip\\
\noindent{\bf Constants for ${\mathbf m}$-homologous cycles on finite orbit surfaces.}\  
Given a $d,n \in \N$, we consider the set 
of saddle connections on  
\[\#^d_{(a+ib)d/n}\T^2 \in \sO_n \subset 
\sF^{sym}_d.\]   
According to the formulas developed earlier 
each saddle connection on 
$\#^d_{(a+ib)d/n}\T^2$ is mapped to a  
horizontal saddle connection on 
some $\#^d_{(j+ik)d/n}\T^2 \in \sO_n$ 
for some $A \in \slz$. This $A$ is 
determined up to a  
parabolic subgroup generated by
$u_1=\left(\begin{smallmatrix} 1&1 \\ 0&1 \end{smallmatrix}\right)$ with orbit   
\[\{ u^m_1 \cdot A\cdot \#^d_{(j+ik)d/n}\T^2: m \in \Z\} 
\cong \Z/w\Z
.\]
The order $w$ of this parabolic group equals the 
{\em width} of the cusp it defines in 
\[\slv{\#^d_{(j+ik)d/n}\T^2,dz} \subset \slz.\] 
Since this width appears as a factor in the 
formula for the Siegel-Veech constants we 
can actually talk about the Siegel-Veech 
constants for saddle connections 
associated to any particular 
surface  (compare equations \ref{rate4} and 
\ref{rate5})
\[ \#^d_{(j+ik)d/n}\T^2 \in \sO_n.\]    
Of course, one cannot really associate 
a particular set of saddle connections, 
i.e. a set of saddle connections with 
distinguished properties depending on  
$\CF_h( \#^d_{(j+ik)d/n}\T^2,dz) $, to 
$ \#^d_{(j+ik)d/n}\T^2 $, 
since these properties are preserved  
under the action of $\{u^n:\ n \in \Z\}$. 
 What we do instead is looking at 
 the Siegel-Veech constant associated to 
 the orbit
\[\sU_{j,k}=\{ u^n_1 \cdot \#^d_{j+ik}\T^2: n \in \N\} 
\subset \partial^{btm}\sC_k. \] 
and divide the result by the cusp-width $|\sU_{j,k}|$. 
In that sense we understand the next theorem. 
Further we note that for a saddle 
connection $s_{a+ib} \in \sS\!\sP_{h,m} $ bounded by 
$\#_{a+ib}\T^2 \in \sO_m$ to the right we 
have $(a,b,d)=m$, 
while $(a+1,b,d)=m$, if   
$\#_{(a+1)+ib}\T^2 \in \sO_m$ bounds 
$s_{a+ib}$ to the left. 
\begin{theorem}
Assume there is an \[\#^d_{t_h+it_v}\T^2 \in 
s_{a+ib} \cap \sO_n, \ \text{ with } s_{a+ib} \in 
\sS\!\sP_{h,m}. \] 
Then the asymptotic quadratic growth rate of the 
$m$-homologous families of saddle connections,  
which are associated to  $\#^d_{t_h+it_v}\T^2$ 
in the sense above,   
contained in $\abas \in \sO_n$ is:   
\begin{equation}
c^d_{a,b,m}(\alpha)=\frac{d}{m \cdot |\sO_n|\cdot |\{t_h\}-\epsilon|^2}.
\end{equation}
Here $\epsilon=1$, if we count saddle 
connections in the orbit of the one 
defined by the horizontal 
segment connecting $\#^d_{t_h+it_v}\T^2$  
with $\#^d_{(a+1)+ib}\T^2 \in \sO_m$  
and $\epsilon=0$ otherwise. 
The quadratic growth rate of all $m$-homologous 
saddle connections on any $\abas \in \sO_n$ 
is:
\begin{equation}
c^d_{m}(\alpha)=\frac{2\cdot d}{m\cdot |\sO_n|}  
\sum_{\#^d_{t_h+it_v}\T^2 \in \sS\!\sP_{h,m} \cap \sO_n} \frac{1}{|\{t_h\}|^2}.
\end{equation}
\end{theorem}
\noindent{\em Note:} Any degeneration of 
  $\#^d_{t_h+it_v}\T^2 \in \sF^{sym}_d$ 
collapsing the two cone points also collapses exactely $d$ 
saddle connections.    
\begin{proof}
Straight forward application of formula \ref{torsadcount}. 
\end{proof}    
Like for cylinders the constants for saddle connections  
become easier if we restrict to torsion points of order $n\neq d$, 
with $(n,d)=1$, because in this case (see equation \ref{baseline}) 
\[\T^2_d(n) \cap \sS\!\sP_{h,m} = 
\T^2_d(n) \cap \partial^{btm}\sC_0.\]
Thus we find for $\#^d_{dj/n}\T^2 \in \T^2_d(n) \cap 
\partial^{btm}\sC_0$ 
\[ c^d_{[dj/n]}(dj/n)=\frac{1}{d\cdot m}\frac{n^2}{\phi(n)\psi(n)} 
\frac{1}{|\{j\}|^2}.\]
Taking the limit $n \rightarrow \infty$ while assuming $(n,d)=1$, 
gives the generic constants: 
\[ c^d_{j,0,\pm}=\frac{\pi^2}{3}\frac{1}{d \cdot m}.\]
\vspace*{2mm}

\noindent {\bf Continuity of Siegel-Veech constants for $d=2$.}\label{convergence}
We discuss convergence of 
Siegel-Veech constants of saddle connections 
connecting the two cone points in case of 
$2$-symmetric differentials. In this case the 
modular fiber is presented by 
\[\sF^{sym}_2 \cong \C/2\Z^2 - \Z^2/2\Z^2.\]
The group $\Z^2/2\Z^2$ parameterizes four 
degenerated surfaces and using the {\em multiplication} 
isomorphism 
\[ m_2: \T^2 \rightarrow \C/2\Z^2 \cong \sF^{sym, c}_2, \]
the three non-vanishing $2$-torsion points on $\T^2$ 
represent one $\slz$ orbit, the other is the point 
$[0] \in \T^2 \cap \frac{1}{2}\Z^2/\Z^2$. 
$\CF_h(\sF^{sym}_2,\omega_2)$ contains $4$ saddle connections 
$s^{\pm}_0,s^{\pm}_1$, where the 
$s^{+}_i$ connect $[i/2,0] \in \T^2 $ with $[i/2, 1/2] \in \T^2$ 
and the $s^{-}_i$ connect  $[i/2,1/2] \in \T^2 $ with $[i/2, 0] \in \T^2$. 
Now for all $n \in \N$ with $(n,2)=1$, 
we have 
\[\sO_n \cap s^{\pm}_1=\{[a/n,b/n] \in \T^2: (a,b,n)=1 \} \cap s^{\pm}_1= 
\emptyset. \]
The Siegel-Veech constant $c^{\pm}_{0}(n)$ associated 
with $[0] \in \T^2$ is not affected by 
this, since $[0]$ is a fixed point of $\slz$. 
By the above the Siegel-Veech constant $c^{\pm}_{1}(n)$ 
for surfaces on the orbit $\sO_n$, $n$ odd, is entirely determined by intersections $\sO_n \cap s^{\pm}_0$. 
In fact, one easily states   
\begin{equation}
c^{\pm}_1(n)=4\frac{n^2}{\varphi(n)\psi(n)}\sum^{(n+1)/2}_{\stackrel{i=1}{(i,2n)=1}}
\frac{1}{i^2}. 
\end{equation}
Now taking the limit over any sequence of {\em odd} numbers
we find
\begin{equation}
c^{\pm}_1=4 \frac{\varphi(2)\psi(2)}{4} \zeta(2)=
3 \cdot \zeta(2).  
\end{equation}
To see this, one can use the Euler product 
\[ \zeta(2)=\sum^{\infty}_{i=1}\frac{1}{i^2}=
\prod_{p \  prime}(1-\frac{1}{p^2})^{-1}=
\prod_{p|q}(1-\frac{1}{p^2})^{-1}\prod_{p \nmid q}(1-\frac{1}{p^2})^{-1},\]
from which one derives (for odd $n$)
\[ \frac{2^2}{\varphi(2)\psi(2)}\frac{n^2}{\varphi(n)\psi(n)}
\sum^{\infty}_{\stackrel{i=1}{ (i,2n)=1}}\frac{1}{i^2} =\zeta(2).\]    
Now take an {\em even} $n$, then 
it follows from  Lemma \ref{torsionformula}, 
Equation \ref{orbitnumber} 
\[|\sO_n \cap s^{\pm}_1|=\varphi(n/2)\ \text{ and }\ 
|\sO_n \cap s^{\pm}_0|=\frac{1}{2}\varphi(n).
 \]
This gives Siegel-Veech constants 
\begin{equation}\label{evencase}
\begin{split}
c^{\pm}_1(n)&=4\frac{n^2}{\varphi(n)\psi(n)}\left[
\frac{1}{2}\sum^{n/2}_{\stackrel{i=1}{(i,n/2)=1}}
\frac{1}{i^2} 
+ \frac{1}{4}
\sum^{n/2}_{\stackrel{i=1}{(i,n)=1}}
\frac{1}{(n/2-i)^2} \right] = \\
&=4\frac{n^2}{\varphi(n)\psi(n)}\begin{cases}   
\frac{1}{2}\sum^{n/2}_{\stackrel{i=1}{(i,n/2)=1}}
\frac{1}{i^2} 
+ \frac{1}{4}
\sum^{n/2}_{\stackrel{i=1}{(i,n)=1}}
\frac{1}{(n/2-i)^2} 
& \text{ if } 4\  | \  n \vspace*{1mm}\\
\frac{1}{2}\sum^{n/2}_{\stackrel{i=1}{(i,n/2)=1}}
\frac{1}{i^2} 
+ \frac{1}{16}
\sum^{n/2}_{\stackrel{i=1}{(i,n)=1}}
\frac{4}{(n/2-i)^2} 
& \text{ if } 4 \nmid n.
\end{cases}
\end{split}
\end{equation}
Under the hypothesis $(i,n)=1$ the difference $n/2-i$ 
can contain any prime-factor, but not the primes 
dividing $n/2$. Thus in case $4|n$ we 
can rearrange and find 
\begin{equation}\label{2evencase}
c^{\pm}_1(n)=
3 \frac{n^2/4}{\varphi(n/2)\psi(n/2)}
\sum^{n/2}_{\stackrel{i=1}{(i,n/2)=1}}
\frac{1}{i^2}  \stackrel{m \rightarrow \infty}
{\longrightarrow} 3 \cdot  \zeta(2)\quad 
\text{ if } 4\  | \  n.
\end{equation}
If $4 \nmid n$ (or $ (n/2,2)=1$ ) 
we have to take into 
account that $2 | (n/2-i)$ and that 
$\frac{n^2}{\varphi(n)\psi(n)}=\frac{4}{3}
\frac{n^2/4}{\varphi(n/2)\psi(n/2)}$
to find
\begin{equation}\label{2oddevencase}
c^{\pm}_1(n)=4\cdot \frac{4}{3}\cdot
\frac{m^2}{\varphi(m)\psi(m)}\left[
\frac{1}{2}\sum^{m}_{\stackrel{i=1}{(i,m)=1}}
\frac{1}{i^2} 
+ \frac{1}{16}
\sum^{m/2}_{\stackrel{i=1}{(i,m)=1}}
\frac{1}{i^2} \right] \stackrel{m \rightarrow \infty}{\longrightarrow}
 3 \cdot \zeta(2) 
\end{equation}
using $m:=n/2$ and Euler products again. 
This establishes the convergence for $2$-symmetric torus 
covers. In the same manner one can show it 
for $d$-symmetric differentials for every $d$. 
We do not write down the details at this place and 
remark that a general convergence theorem, 
i.e. one which just depends on the Siegel-Veech 
formula and Ratner's Theorem, would be 
the desired and expected result. \label{convergenceend}
\medskip\\
\noindent{\bf Comments.} The example treated in this paper as well as the various asymptotic formulas show that the asymptotic constants 
are very sensitive invariants of the $\slz$ orbit of a 
particular translation surface in $\sF^{sym}_d$. The previous 
discussion of saddle connections shows, that we can 
characterize cusps of a particular $\slz$-orbit using 
particular types of saddle connections.  The asymptotic quadratic constants 
are good enough to distinguish (finite) $\slz$-orbits, moreover 
if one adds the horizontal cylinder decomposition of the 
modular fiber, the quadratic constants are good enough to 
distinguish cusps of a particular $\slz$ orbit.  
\medskip\\
As a matter of fact the quadratic growth rates of saddle conntections 
depend only on the saddle connections in the horizontal foliation of 
$\sF^{sym}_d$. The generic constants however, obtained 
by approximation along the horizontal foliation  
are also the constants for saddle connections which will never be horizontal under application 
of an $A \in \slz$.


\vspace*{10mm}

\end{document}